\documentclass[12pt]{amsart}

\usepackage{amssymb}
\usepackage{amsmath}
\usepackage{mathrsfs}

\newcommand{\Lg}{\mbox{$\mathfrak g$}}
\newcommand{\Lgr}{\mbox{$\mathfrak g_{\mathbf R}$}}
\newcommand{\Lkr}{\mbox{$\mathfrak k_{\mathbf R}$}}
\newcommand{\Lpr}{\mbox{$\mathfrak p_{\mathbf R}$}}
\newcommand{\Lhg}{\mbox{$\hat{\mathfrak g}$}}
\newcommand{\Lhgr}{\mbox{$\hat{\mathfrak g}_{\mathbf R}$}}
\newcommand{\Lh}{\mbox{$\mathfrak h$}}
\newcommand{\Lk}{\mbox{$\mathfrak k$}}
\newcommand{\Lp}{\mbox{$\mathfrak p$}}

\newcommand{\Lm}{\mbox{$\mathfrak m$}}

\newcommand{\Lu}{\mbox{$\mathfrak u$}}
\newcommand{\Lz}{\mbox{$\mathfrak z$}}

\newcommand{\Pf}{{\em Proof}. }
\newcommand{\EPf}{\hfill$\square$}

\newcommand{\R}{\mbox{$\mathbf R$}}
\newcommand{\C}{\mbox{$\mathbf C$}}

\newcommand{\inn}[2]{\mbox{$\mathcal{h} #1,#2 \mathcal{i}$}}
\newcommand{\ad}[1]{\mbox{$\mbox{ad}_{#1}$}}
\newcommand{\Ad}[1]{\mbox{$\mbox{Ad}_{#1}$}}

\newcommand{\Spin}[1]{\mbox{$\mathbf{Spin}(#1)$}}

\newcommand{\Gr}{\mbox{$G_{\mathbf R}$}}
\newcommand{\hGr}{\mbox{$\hat{G}_{\mathbf R}$}}
\newcommand{\Vr}{\mbox{$V_{\mathbf R}$}}
\newcommand{\Kr}{\mbox{$K_{\mathbf R}$}}
\newcommand{\tth}{\mbox{$\tilde\theta$}}
\newcommand{\ts}{\mbox{$\tilde\sigma$}}
\newcommand{\tom}{\mbox{$\tilde\omega$}}
\newcommand{\tmu}{\mbox{$\tilde\mu$}}
\newcommand{\cs}{\mbox{$c^{\tilde\sigma}$}}
\newcommand{\ct}{\mbox{$c^{\tilde\theta}$}}
\newcommand{\cmt}{\mbox{$c^{-\tilde\theta}$}}
\newcommand{\cst}{\mbox{$\cs\cap\ct$}}
\newcommand{\csmt}{\mbox{$\cs\cap\cmt$}}
\newcommand{\cat}{\mbox{$/\!/ $}}

\newcommand{\dimr}{\mbox{$\dim_{\mathbf R}$}}

\newtheorem{thm}{Theorem}

\newtheorem{cor}[thm]{Corollary}
\newtheorem{prop}[thm]{Proposition}
\newtheorem{lem}[thm]{Lemma}
\newtheorem{eg}[thm]{Example}
\newtheorem{conj}[thm]{Conjecture}

\newtheorem{claim}{\sc Claim}

\begin{document}

\title{Polar orthogonal representations of\\
real reductive algebraic groups}
\author{Laura Geatti and Claudio Gorodski}
\address{Dipartimento di Matematica, Universit\`a di Roma 2 Tor Vergata,
         via della Ricerca Scientifica, 00133 Roma, Italy}
\email{geatti@mat.uniroma2.it}
\address{Instituto de Matem\' atica e Estat\'\i stica,
         Universidade de S\~ ao Paulo,
         Rua do Mat\~ao, 1010,
         S\~ ao Paulo, SP 05508-090,
         Brazil}
\email{gorodski@ime.usp.br}
%\thanks{The second author was supported in part by 
%FAPESP and CNPq.}

\date{\today}

\begin{abstract}
We prove that a polar orthogonal representation of a real reductive
algebraic group has the same closed orbits as the 
isotropy representation of a pseudo-Riemannian symmetric space. 
We also develop a partial structural theory of polar 
orthogonal representations of real reductive
algebraic groups which slightly generalizes some results
of the structural 
theory of real reductive Lie algebras.
\end{abstract}

\maketitle

\section{Introduction}

A representation of a complex reductive algebraic
group~$G$ on a finite-dimensional complex vector space~$V$ is called
\emph{polar} if there exists a subspace $c\subset V$ consisting
of semisimple elements such that $\dim c = \dim V\cat G$ (the categorical
quotient), and for a dense subset of $c$, 
the tangent spaces to the orbits are 
parallel~\cite{DK}; then it turns out that every closed orbit 
of $G$ meets $c$ (Prop.~2.2, \emph{loc.\ cit.\/}). The class
of polar representations was introduced and studied 
by Dadok and Kac in~\cite{DK}, 
and it is very important in invariant theory
because it includes the adjoint actions, the representations
associated to symmetric spaces studied by Kostant and Rallis~\cite{KR}  
as well as, more generally, the representations associated to 
automorphisms of finite order ($\theta$-groups) introduced by 
Vinberg~\cite{V} (see also~\cite{Kac2}). At present, there is 
no complete classification of polar representations although
the paper~\cite{DK} contains very important partial results. 

A complex (resp.~real) representation admitting a complex-valued
(resp.~real-valued) invariant non-degenerate
symmetric bilinear form is called \emph{orthogonal}. 
It is well known that a complex orthogonal representation 
admits a real form
invariant under a maximal compact subgroup. Consider in particular
the complex polar orthogonal representations and the class of compact 
real forms they originate. 
Since the complex reductive algebraic 
groups are exactly the complexifications of the compact Lie groups,
one can equivalently
define directly the concept of a real polar representation
of a compact Lie group in the differential-geometric setting
(as in e.g.~\cite{PT}) and obtain the same class.
Note that orbits of polar representations of compact 
Lie groups are very important in submanifold geometry and Morse 
theory~\cite{BS,C,Sz,PT,DO,GTh3}. Now, such representations were 
classified by Dadok in~\cite{D}, and the following very nice 
characterization was deduced:
\emph{A polar representation of a compact Lie group has the same
orbits as the isotropy representation of a Riemannian
symmetric space.}

The purpose of this paper is to study 
\emph{non-compact} real forms of complex polar orthogonal 
representations.
Equivalently, we define a representation of 
a real reductive algebraic group (in the sense
of~\cite[\S1]{BH}) to be \emph{polar} 
if and only if its algebraic complexification is polar.
In section~\ref{sec:classif} we prove the following theorem.

\begin{thm}\label{classif}
A polar orthogonal representation of a connected
real reductive 
algebraic group has the same closed orbits as the isotropy
representation of a pseudo-Riemannian symmetric space.
\end{thm}

In section~\ref{sec:isoparam}, 
we discuss some aspects of the submanifold 
geometry of the closed orbits of the polar orthogonal representations
of the real reductive algebraic groups in that we relate 
them to a notion of pseudo-Riemannian isoparametric submanifold
of a pseudo-Euclidean space (compare~\cite{Hahn,Magid}). 
%(Theorem~\ref{isoparam}). 

Finally in section~\ref{sec:structure}, 
independently of classification results, we develop
a partial structural theory of polar 
orthogonal representations of real reductive 
algebraic groups that generalizes some results of the structural 
theory of real reductive Lie algebras. In this regard,
we propose to replace adjoint actions by  
polar orthogonal ones. The results we prove are slight 
generalizations
of well known results for the adjoint actions, but we believe
our proofs are more geometric. 
In particular, we show that 
a polar orthogonal representation of a real reductive 
algebraic group admits finitely many pairwise inequivalent 
so called Cartan subspaces in standard position 
with respect to a compact real form
such that the union of those
subspaces meets all the closed orbits and always orthogonally
(Theorem~\ref{finiteness} and Corollary~\ref{cor:intersection}). 
We also construct the so called
Cayley transformations that relate different equivalence classes
of Cartan subspaces (\S\ref{cayley}), 
and use those to show that  
the equivalence classes of Cartan subspaces in the 
two extremal positions with respect to 
the compact real form are unique (Corollary~\ref{extremal}). 

Unless explicit mention to
the Zariski topology is made, we use throughout 
the classical topology. 
We always use
lowercase gothic letters to denote Lie algebras.
For a given homomorphism of groups, we denote 
the induced homomorphism on the Lie algebra level by the same 
letter whenever the context is clear.
Sometimes it is useful to call 
a representation \emph{orthogonalizable} if it admits an 
invariant non-degenerate symmetric bilinear form 
but we do not want to fix such a form. 

Part of this work was completed while the first author was visiting
University of S\~ao Paulo for which she would like to thank CNPq
and CCInt-USP for their generous support, and part of it 
was completed while the second author was visiting Universit\`a
di Roma 2 ``Tor Vergata'' for which he would like to thank
INdAM for its generous support. The second author also wishes
to thank Ralph Bremigan and Martin Magid 
for very useful discussions.

\section{Preliminaries}\label{prelim}

Let $G$ be a connected complex reductive algebraic group. 
Let $\tau:G\to GL(V)$ be a complex representation. 
A vector $v\in V$ is called \emph{semisimple} if the orbit $Gv$ is closed. 
Not every orbit of $G$ in $V$ is closed, 
but the closure of any orbit contains a unique 
closed orbit. An element is called \emph{regular}
if it is semisimple and $\dim Gv\geq \dim Gx$ for all
semisimple $x\in V$. The representation $\tau$ is called 
\emph{stable} or is said to admit \emph{generically closed
orbits} if there exists an open and dense 
subset of $V$ consisting of closed orbits. 
An orthogonalizable representation is necessarily 
stable (see~\cite[Cor.~5.9]{S} or~\cite{Lu2,Lu3}). 

Let $\C[V]$ be the polynomial algebra of $V$, and let 
$\C[V]^G$ be the algebra of $G$-invariant polynomials.
It does not contain nilpotents, and is finitely generated
by a theorem of Hilbert, so it is the coordinate ring of 
an affine algebraic variety denoted by $V\cat G$ and called
the \emph{categorical quotient} of $V$ by $G$. 
The embedding $\C[V]^G\to\C[V]$ 
induces a surjective morphism of affine algebraic varieties
$\pi:V\to V\cat G$. Every fiber of $\pi$ contains a unique
closed orbit. It follows that $V\cat G$ can be seen as 
the parameter set of closed $G$-orbits in $V$, and 
then $\pi(v)$ represents the unique closed orbit 
in the closure of $Gv$~\cite[\S4]{PV}. 

For semisimple $v\in V$, set
\[ c_v = \{\, x\in V \;|\; \Lg\cdot x\subset\Lg\cdot v\,\}. \]
Then $c_v$ consists entirely of semisimple 
elements, and the isotropy subalgebras satisfy
$\Lg_x\supset\Lg_v$ for~$x\in c_v$~\cite[Lem.~2.1]{DK}. 
The representation $\tau$ 
is called \emph{polar} if a semisimple $v$ can be chosen
so that $\dim c_v = \dim V\cat G$. In this case, $c_v$ is called
a \emph{Cartan subspace}. The Cartan subspaces of a polar
representation are all $G$-conjugate~\cite[Thm.~2.3]{DK}.

The group $G$ can be seen simply as the complexification of a 
compact connected Lie group $U$; compare~\cite[\S5]{S} 
or~\cite[Rmk.~3.4]{BH}.
Then $U$ is a maximal compact (necessarily connected) subgroup of $G$, 
and every maximal compact subgroup of $G$ is $G$-conjugate to $U$. 
It is easy to see that a representation $\tau$ is orthogonalizable
if and only if it admits a real form 
$\tau_u:U\to GL(W)$~\cite[Prop.~5.7]{S}. The group $U$ must be the fixed
point group $G^\theta$ of a unique anti-holomorphic 
involutive automorphism
$\theta$ of $G$, which is called a \emph{Cartan involution}
of $G$.
Also, the subspace $W$ is the fixed point set $V^{\tilde\theta}$ 
of a conjugate-linear involutive automorphism  
$\tilde\theta$ of $V$, the equation
$\tilde\theta(g\cdot v) = \theta(g)\cdot\tilde\theta(v)$
holds for $g\in G$ and $v\in V$, and an invariant form
$\inn{\cdot}{\cdot}$ can be chosen on $V$
so that it is real-valued on $V^{\tilde\theta}$. 
 
More generally, we consider real forms of  
$\tau:G\to GL(V)$ given by a pair $(\sigma,\ts)$ where 
$\sigma$ is an anti-holomorphic 
involution of $G$ and $\ts$ is a real structure on $V$ satisfying 
$\ts(g\cdot v)=\sigma(g)\cdot\ts(v)$ for $g\in G$, $v\in V$.
The fixed point subgroup $G^\sigma$ is a 
(not necessarily connected)
real reductive algebraic group, and $\tau$ of course restricts to
a representation of $G^\sigma\to GL(V^{\tilde\sigma})$, where
$V^{\tilde\sigma}$ is the fixed point set of $\ts$ in $V$. 
We say that two real forms $(\sigma,\tilde\sigma)$ and 
$(\sigma',\tilde\sigma')$ \emph{commute}
if they commute componentwise. 
If $\tau$ is 
orthogonal with respect to $\inn{\cdot}{\cdot}$
and a real form $(\sigma,\ts)$ is given,
then $\inn{\cdot}{\cdot}$ is said to be 
\emph{defined over $\R$ with respect to
$\tilde\sigma$} and $(\sigma,\ts)$ is called an 
\emph{orthogonal real form} if $\inn{\cdot}{\cdot}$ is real-valued
on $V^{\tilde\sigma}$. Note that the latter condition is equivalent to
having 
\[ \inn{\tilde\sigma x}{\tilde\sigma y}=\overline{\inn xy} \]
for $x$, $y\in V$. 
A \emph{Cartan pair} of $\tau$ is an orthogonal 
real form $(\theta,\tilde\theta)$ such that $\theta$ is a Cartan  
involution of $G$ and $\inn{\cdot}{\cdot}$ 
is real-valued and negative-definite on
$V^{\tilde\theta}$. Note that $(\theta,\tth)$ is a 
Cartan pair of $\tau$ with respect to $\inn{\cdot}{\cdot}$ if and 
only if $(\theta,-\tth)$ is a 
Cartan pair of $\tau$ with respect to $-\inn{\cdot}{\cdot}$.
The following result is essentially
proved in~\cite[7.4]{Br}, but we find it convenient to include a 
proof here because we will need to refer to some of its techniques. 

\begin{prop}[Bremigan]\label{cartan-pairs}
Let $\tau:G\to O(V,\inn{\cdot}{\cdot})$ be an orthogonal representation, and 
suppose that $(\sigma,\ts)$ is an orthogonal real form.
Then there exists a Cartan pair $(\theta,\tth)$ 
which commutes with $(\sigma,\ts)$. 
\end{prop}

\Pf It is well known that 
there exists a Cartan involution $\theta$ of $G$ such that
$\theta\sigma=\sigma\theta$. Let $U=G^\theta$ be the associated
maximal compact subgroup of $G$. Consider the realification $V^r$ of $V$,
and denote the invariant complex structure on $V^r$ by $J$ so that 
$V=(V^r,J)$. Note that 
\[ \tilde\sigma\tau(g)\tilde\sigma^{-1}=\tau(\sigma(g)),
\quad J\tau(g)J^{-1}=\tau(g)\quad\mbox{and}\quad J\tilde\sigma J^{-1}
=-\tilde\sigma \]
for $g\in G$. Let $G^*$ be the subgroup of $GL(V^r)$ generated
by $\tau(G)$, $\tilde\sigma$ and $J$. Then $G^*$ contains $\tau(G)$ 
as a normal subgroup of finite index. Due to $\theta\sigma=\sigma\theta$,
we have also that $\tilde\sigma$ normalizes $\tau(U)$. Let 
$U^*$ be the subgroup of $G^*$ generated by $\tau(U)$, $\tilde\sigma$ 
and $J$. 
Then $U^*$ is a compact subgroup of $G^*$, so we can find an 
$U^*$-invariant 
positive-definite real inner product on $V^r$ which we denote by 
``$\cdot$''. Set 
\[ (x,y)=x\cdot y + i (x\cdot Jy) \]
for $x$, $y\in V^r$. Then it is easily checked that
$(\cdot,\cdot)$ is an $U$-invariant positive-definite 
Hermitian form on $(V^r,J)=V$ which is real-valued on $V^{\tilde\sigma}$.
In particular, $i\Lu$ acts on $V$ by Hermitian endomorphisms.  
Next, define a conjugate-linear automorphism $\tth$ of $V$ by 
setting 
\begin{equation}\label{hermitian-symmetric}
 (x,\tilde\theta y)=-\inn xy 
\end{equation}
for $x$, $y\in V$. Then 
\begin{equation}\label{3}
 (x,\tth(gy))=-\inn{x}{gy}=-\inn{g^{-1}x}y=(g^{-1}x,\tilde\theta y)
=(x,\theta(g)\tilde\theta(y)), 
\end{equation}
so $\tilde\theta\tau(g)=\tau(\theta(g))\tth$
for $g\in G$. Moreover
\begin{equation}\label{4}
 (x,\tth\ts y)=-\inn x{\ts y}=-\overline{\inn{\ts x}y}
=\overline{(\ts x,\tth y)}=(x,\ts\tth y) 
\end{equation}
implying that $\tth\ts=\ts\tth$. 
We also have that 
\begin{equation*}
\begin{split}
(\tth^2 x,y)&=\overline{(y,\tth^2 x)}=
-\overline{\inn y{\tth x}}=-\overline{\inn{\tth x}y}\\
&\quad=\overline{(\tth x,\tth y)}
=(\tth y,\tth x)=\cdots\\
&\quad=(x,\tth^2y)
\end{split}
\end{equation*}
for $x$, $y\in V$.
It follows that $\tilde\theta^2:V\to V$ is a $G^*$-equivariant 
$\C$-linear Hermitian automorphism. Hence there exists a 
$\inn{\cdot}{\cdot}$- and $(\cdot,\cdot)$-orthogonal 
$G^*$-invariant decomposition 
$V=\oplus_j V_j$ such that $\tth^2|_{V_j}=\lambda_j\,\mathrm{id}_{V_j}$
where $\lambda_j\in\R\setminus\{0\}$ and the $\lambda_j$'s are
pairwise distinct. 

Note that $\lambda_j(x,x)=(\tth x,\tth x)>0$
if $x\in V_j\setminus\{0\}$, so we also have $\lambda_j>0$. 
If we change $(\cdot,\cdot)$ by a factor of
$\lambda_j^{1/2}$ on $V_j\times V_j$, as we do,
$\tth$ is changed by a factor of $\lambda_j^{-1/2}$ on $V_j$, and then 
the resulting $\tth$ satisfies $\tth^2=\mathrm{id}_V$. 
Note that equations~(\ref{3}) and~(\ref{4}) are unchanged.
Now $(\theta,\tth)$ is a real form of $(G,V)$ commuting
with $(\sigma,\ts)$. Further,
\[ \inn{\tth x}{\tth y}=-(\tth x,\tth^2y)=-(\tth x,y)=-\overline{(y,\tth x)}
=\overline{\inn yx}=\overline{\inn xy} \]
for $x$, $y\in V$ and 
\[ \inn xx=-(x,\tth x)=-(x,x)<0 \]
for $x\in V^{\tilde\theta}\setminus\{0\}$. This completes
the proof. \EPf

\begin{prop}\label{cartan-pairs-uniqueness}
Let $\tau:G\to O(V,\inn{\cdot}{\cdot})$ be an orthogonal representation, and 
suppose that $\theta$ is a Cartan involution of $G$. Then there can be
at most one real structure $\tth$ on $V$ such that 
$(\theta,\tth)$ is a Cartan pair of $\tau$. 
\end{prop}

\Pf Suppose
that $(\theta,\tth)$ and $(\theta,\tth')$ are two 
Cartan pairs of $\tau$.
Define 
\[ h(x,y) = -\inn{x}{\tth y}\quad\mbox{and}\quad h'(x,y)=-\inn{x}{\tth'y} \]
for $x$, $y\in V$. It is easy to see that $h$ and $h'$ are
two $U$-invariant
positive-definite Hermitian forms.
%\ft{These will be real-valued 
%on~$V^{\tilde\sigma}$ if $\tth$, $\tth'$ commute with a real
%structure $\ts$ on $V$.}
Diagonalizing $h'$ with respect to $h$, we get a 
$U$-invariant, $h$-orthogonal splitting $V=\oplus_j V_j$ such that 
$\tth'=\lambda_j\tth$ on $V_j$, where 
$\lambda_j>0$ and the $\lambda_j$'s are pairwise distinct. 
Using $(\tth')^2=\tth^2=1$, we finally see 
that $\tth'=\tth$. \EPf

\begin{cor}\label{conjugate-cartan-pairs}
Let $\tau:G\to O(V,\inn{\cdot}{\cdot})$ be an orthogonal representation.
Then any two Cartan pairs of $\tau$ are $G$-conjugate; 
moreover, if the underlying Cartan involutions 
commute with a real form $\sigma$ of $G$, 
then the Cartan pairs are $(G^\sigma)^\circ$-conjugate.
\end{cor}

\Pf Let $(\theta_1,\tilde\theta_1)$ and 
$(\theta_2,\tilde\theta_2)$ be two Cartan pairs of $\tau$. 
It is known that there exists $g\in G$ such that 
$\theta_2=\mathrm{Inn}_g\theta_1\mathrm{Inn}_g^{-1}$. Of course, 
$(\mathrm{Inn}_g\theta_1\mathrm{Inn}_g^{-1},g\tth_1g^{-1})$ is also a Cartan pair.
Proposition~\ref{cartan-pairs-uniqueness} implies that 
$\tth_2=g\tth_1g^{-1}$. Further,
if both $\theta_1$ and $\theta_2$ commute with $\sigma$,
it is known that $g$ can be taken in the identity
component of $G^\sigma$. \EPf

\section{The classification}\label{sec:classif}

Let $\hGr/\Gr$ be a 
semisimple pseudo-Riemannian symmetric space. 
Here $\hGr$ is a connected real semisimple Lie group, 
$\hat\tau$ is a non-trivial
involutive automorphism of $\hGr$ and $\Gr$ is 
an open subgroup of the fixed point group
of $\hat\tau$. The automorphism
$\hat\tau$ induces an automorphism of the Lie algebra $\Lhgr$ of $\hGr$
which we denote by the same letter. Let $\Lhgr=\Lgr+\Vr$ be the 
decomposition into $\pm1$-eigenspaces of $\hat\tau$. Of course,
$\Lgr$ is the Lie algebra of $\Gr$. The restriction of the 
Killing form of $\Lhgr$ to $\Vr\times\Vr$ is $\mathrm{Ad}_{G_{\mathbf R}}$-invariant
and non-degenerate, so it induces a $\hGr$-invariant pseudo-Riemannian metric
on $\hGr/\Gr$. The adjoint action of $\Gr$ on $\Vr$ is equivalent to 
the isotropy representation of $\hGr/\Gr$ at the base-point.

Next, extend $\hat\tau$ complex-linearly to an automorphism
of the complexification $\Lhg=(\Lhgr)^c$ denoted by the 
same letter and 
consider the corresponding decomposition
$\Lhg=\Lg+V$ into $\pm1$-eigenspaces. Let $\hat G$ be the simply-connected
complex Lie group with Lie algebra $\Lhg$, view $\hat\tau$
as an involution of $\hat G$, and let $G$ be the 
fixed point group of $\hat\tau$ in $\hat G$. Note that $G$ is connected.  
The adjoint action of $G$ on $V$ 
is a complex polar action whose Cartan subspaces coincide with  
the maximal Abelian subspaces of $V$ consisting of 
semisimple elements
(indeed, this is a $\theta$-group (see~\cite[Introd.]{DK}
or~\cite[8.5, 8.6]{PV}; no relation here to the aforementioned
Cartan involution $\theta$). Further, it is an 
orthogonal action with respect to the restriction of the 
Killing form of $\Lhg$ to~$V$. 
By passing from $\hGr$ to a finite covering if necessary,
we may assume that $\hGr$ embeds into $\hat G$ and $\Gr$ embeds
into $G$, so we can view 
the adjoint action of $\Gr$ on $\Vr$ as an orthogonal real form
of the adjoint action of $G$ on $V$. 
%Recall also that the closed orbits of $G$ on $V$ 
%(resp.~$\Gr$ on $\Vr$) are precisely those through semisimple 
%elements, where $X\in\Ll$ is called semisimple if and only 
%if $\mathrm{ad}_X$ is a semisimple endomorphism of $\Ll$. 
%Therefore this terminology is consistent with that introduced 
%in section~\ref{prelim}. 
We deduce 
that the isotropy representation 
of a pseudo-Riemannian symmetric space is a polar representation.
In this section, we prove Theorem~\ref{classif} which is 
essentially a converse to this result.  

Before giving the proof of Theorem~\ref{classif}, we 
prove four lemmas. 
In the remaining of this section, let $G$ be a complex reductive
algebraic group defined over~$\R$ and denote by 
$\Gr$ the identity component of its real points. 

\begin{lem}\label{complex-polar-reducible}
Let $\tau:G\to GL(V)$ be a polar representation, where 
$V=V_1\oplus V_2$ is a $G$-invariant decomposition. Assume that  
the induced representations $\tau_i:G\to GL(V_i)$ are stable. 
Then:
\begin{enumerate}
\item $\tau_i$ is polar, $i=1$, $2$. 
\item Every Cartan subspace of $\tau$ is of the form 
$c=c_1\oplus c_2$, where $c_i$ is a Cartan subspace of $\tau_i$,
$i=1$, $2$.
\item The closed orbits of $G$ on $V_2$ coincide with those 
of $G_{v_1}$, where $v_1$ is any semisimple vector of $V_1$.
In particular, $V_1$ and $V_2$ are inequivalent representations.
\item Fix a Cartan subspace $c=c_1\oplus c_2$, let 
$\Lh_1$ (resp.~$\Lh_2$) denote the centralizer of $c_2$ (resp.~$c_1$)
in $\Lg$, and denote by $H_i$ the connected subgroup of $G$ 
corresponding to $\Lh_i$. Then the closed orbits of $\tau$
coincide with those of $\hat\tau:H_1\times H_2\to GL(V_1\oplus V_2)$,
where $\hat\tau(g_1,g_2)(v_1+v_2)=\tau_1(g_1)v_1+\tau_2(g_2)v_2$. 
\end{enumerate}
\end{lem}

\Pf Parts~(a) and~(b) are Prop.~2.14 in~\cite{DK}, and~(c)
is Cor.~2.15 of that paper. Let us prove~(d). Select a 
regular element $v_1+v_2\in c_1\oplus c_2$ for $\tau$. Then 
$\Lh_1=\Lg_{v_2}$, $\Lh_2=\Lg_{v_1}$ and (c) implies
that $\Lh_1\cdot v_1=\Lg\cdot v_1$ and
$\Lh_2\cdot v_2=\Lg\cdot v_2$
(note that $v_i$ is semisimple for $\tau_i$ since $c_i$ is a 
Cartan subspace). This implies that $\Lg=\Lh_1+\Lh_2$, so
$G=H_1\cdot H_2=H_2\cdot H_1$ by connectedness of $G$. For any
$u_1+u_2\in c_1\oplus c_2$, we now have that 
$G(u_1+u_2)\subset(H_1\times H_2)(u_1+u_2)$. 
Since $\Lg=\Lh_1+\Lh_2$, and $G(u_1+u_2)$, $H_1 u_1\times
H_2 u_2$ are closed and connected, it follows that the two orbits
coincide. \EPf

\begin{lem}\label{real-polar-reducible}
Let $\rho:\Gr\to GL(\Vr)$ be a polar representation, 
where $\Vr=(\Vr)_1\oplus(\Vr)_2$ is a 
$\Gr$-invariant decomposition. Assume that  
the induced representations $\rho_i:\Gr\to GL((\Vr)_i)$ are orthogonalizable. 
Then there exist closed connected subgroups $H_i'$ of $\Gr$, $i=1$, $2$,
such that the restricted representations 
$\rho_i|_{H_i'}:H_i'\to GL((\Vr)_i)$ are polar 
and the closed orbits of $\rho$ coincide with those
of $\hat\rho:H_1'\times H_2'\to GL((\Vr)_1\oplus(\Vr)_2)$,
where $\hat\rho(g_1,g_2)(v_1+v_2)=\rho_1(g_1)v_1+\rho_2(g_2)v_2$. 
\end{lem}

\Pf The complexification $\tau=\rho^c:G\to GL(V)$ is polar
and each $\tau_i=\rho_i^c$ is orthogonalizable, hence stable. 
By Lemma~\ref{complex-polar-reducible}, 
the closed orbits of $\tau$
coincide with those of $\hat\tau:H_1\times H_2\to GL(V_1\oplus V_2)$,
where $V_i=(\Vr)_i^c$, the group
$H_1$ (resp.~$H_2$) is the connected 
centralizer of $c_2$ (resp.~$c_1$) in $G$, 
and $c=c_1\oplus c_2\subset V_1\oplus V_2$ 
is a Cartan subspace of $\tau$. 
As usual, suppose that $\rho$ is defined by $(\sigma,\ts)$.
Now $c$ can be taken to be 
$\tilde\sigma$-stable due to Lemma~\ref{ts-stable-cartan-subspaces} below. 
In this case, $H_i$ is $\sigma$-stable;
set $H_i'$ to be subgroup of $\Gr$ given by the identity
component of $(H_i)^\sigma$. It is clear that 
the groups $H_i'$ have the desired properties. \EPf

\medskip

Given a representation $\rho:\Gr\to GL(\Vr)$, denote by 
$\rho^*:\Gr\to GL(\Vr^*)$ the dual representation. 
Note that $\rho\oplus\rho^*:\Gr\to GL(\Vr\oplus\Vr^*)$ is always
orthogonal with respect to 
\begin{equation}\label{inner}
 \inn{(v_1,v_1^*)}{(v_2,v_2^*)}=v_1^*(v_2)+v_2^*(v_1). 
\end{equation}
The proof of the following lemma is simple and we omit it. 

\begin{lem}\label{decomp}
Let $\rho:\Gr\to GL(\Vr)$ be orthogonalizable. Then there exists an 
irreducible
decomposition 
\[ \Vr=(\Vr)_{1}\oplus\cdots\oplus(\Vr)_r\oplus
(\Vr)_{r+1}\oplus(\Vr)_{r+1}^*\oplus\cdots\oplus(\Vr)_s\oplus(\Vr)_s^*, \]
where $(\Vr)_{1},\ldots,(\Vr)_r$ are orthogonalizable and 
$(\Vr)_{r+1},\ldots,(\Vr)_s$ are not orthogonalizable. 
\end{lem}

The following lemma will be used to show that certain
polar representations have the same closed orbits as
a the isotropy representation of a symmetric space. 

\begin{lem}\label{same-closed-orbits}
Suppose $\tau:G\to GL(V)$ is a polar orthogonalizable 
representation, $U$ is a maximal compact subgroup of $G$
and $\tau_u:U\to GL(W)$ is a real form. Suppose also
that $U'$ is a connected 
closed subgroup of $U$ and $G'\subset G$ is the 
complexification of $U'$. If $\tau_u|_{U'}$ has the 
same orbits in $W$ as $\tau_u$, then $\tau|_{G'}$ has the same closed
orbits in $V$ as $\tau$. If, in addition, 
$\rho:\Gr\to GL(\Vr)$ is a real form of $\tau$ and
$G'_{\mathbf R}\subset\Gr$ is a connected real form of $G'$, then
$\rho|_{G'_{\mathbf R}}$ has the same closed
orbits in~$\Vr$ as~$\rho$.
\end{lem}

\Pf The assertion about $\rho$ immediately follows
from that about $\tau$ and the facts that 
$\Gr v$ is closed if and only if $Gv$ is closed~\cite{Bi}
and $\dim_{\mathbf R}\Gr v = \dim Gv$ for $v\in\Vr$. 
Let us prove the assertion about $\tau$. 
We first claim that if $v\in V$ and $Gv$ is closed, then
$G'v=Gv$. Of course, we already have that $G'v\subset Gv$. 
In the case in which $v\in W$, we have that
both $Gv$ and $G'v$ 
are connected, closed and have dimension equal 
to $\dimr Uv=\dimr U'v$, so the result follows. 
In the general case, fix a $U$-invariant positive-definite 
Hermitian form $(\cdot,\cdot)$ and 
choose $v_1\in Gv$ of minimal length~\cite[p.508]{DK}.
Of course, $Gv_1=Gv$ and $v_1$ is also of minimal length 
in $G'v_1$. It follows that $G'v_1$ is also 
closed~\cite[Thm.~1.1]{DK} and $G_{v_1}$, $G'_{v_1}$
are $\theta$-stable, where $\theta$ is the  
Cartan involution of $G$ associated to~$U$~\cite[Prop.~1.3]{DK}. 
Let $L=(G_{v_1})^\theta$ and $L'=(G'_{v_1})^\theta$.
Now we can choose $w\in W$ such that $U_w=L$ by the 
same argument as in~\cite[Prop.~5.8]{S}, and it easily follows
that $U'_w=L'$. We have established that $U_w$ (resp.~$U'_w$)
is a real form of $G_{v_1}$ (resp.~$G'_{v_1}$). 
Therefore
\begin{eqnarray*}
    \dim Gv_1 & = & \dim G-\dim G_{v_1} \\
   & = & \dimr U - \dimr U_w \\
 & = & \dimr Uw \\
 & = & \dimr U'w \\
 & = & \dimr U'- \dimr U'_w \\
 & = & \dim G' - \dim G'_{v_1} \\
 & = & \dim G'v_1,
\end{eqnarray*} 
which implies that $G'v_1=Gv_1$. Since $Gv_1=Gv$, we also
have $G'v=Gv$, proving the claim.

Let $c\subset V$ be a Cartan subspace of $\tau$. 
In view of the claim proved above, $c$ consists 
of semisimple elements of $\tau|_{G'}$. Also, 
$\dim c= \dim V\cat G=\dim V\cat G'$, where the last
equality follows from the fact that $\tau_u$ and $\tau_u|_{U'}$
have the same co-homogeneity in $W$. By~\cite[Prop.~2.2]{DK},
every closed $G'$-orbit meets $c$, from which 
it follows that $\tau|_{G'}$ has the same closed orbits in $V$
as $\tau$. \EPf

\medskip

In order to prove Theorem~\ref{classif}, we will use
the explicit lists of polar representations of compact Lie groups
that have been obtained in~\cite{EH1} (irreducible case) 
and~\cite{Ber,Ber2}
(reducible case); see also~\cite{GTh1} (both cases). 
For brevity, an isotropy representation of a semisimple
symmetric space will be called an \emph{s-representation}.  
Let $\rho:\Gr\to GL(\Vr)$ be a polar orthogonal representation. 
Let $\tau=\rho^c:G\to GL(V)$, and suppose that $\rho$ is given by
$(\sigma,\ts)$ so that $\Gr$ is the identity component 
of $G^\sigma$. Let $(\theta,\tth)$ be a Cartan pair as 
in Proposition~\ref{cartan-pairs}, $U=G^\theta$, $W=V^{\tilde\theta}$,
and $\tau_u:U\to GL(W)$ the associated real form. Then $\tau_u$ is a polar
representation of a compact Lie group. By Dadok's theorem 
quoted in the introduction and 
the results in~\cite{EH1,Ber2}, 
$\tau_u$ is either a Riemannian s-representation or one of 
the exceptions
listed in those papers. We need the following 
fundamental lemma. 

\begin{lem}\label{irred-s-repr}
If $\tau_u$ is an irreducible Riemannian s-representation, 
then $\rho$ is a pseudo-Riemannian s-representation. 
\end{lem}

\Pf By assumption, $\hat\Lu=\Lu+W$ admits a 
real Lie algebra structure
extending that of $\Lu$ such that~\cite[p.182]{HZ}
\[ [X,w] = \tau_u(X)w\quad\mbox{and}\quad 
\langle X,[w,w']\rangle_{\mathfrak u} = \inn{\tau_u(X)w}{w'} \]
for $X\in\Lu$ and $w$, $w'\in W$, where 
\[ \langle X,Y\rangle_{\mathfrak u} = 
\mathrm{trace}_{\hat{\mathfrak u}}(\ad X\ad Y) \]
for $X$, $Y\in\Lu$ and $\ad X(Z)=[X,Z]$ for $X\in\Lu$ and $Z\in\hat\Lu$.
Denote the Killing form of $\hat\Lu$ by $\beta$;
note that it is nondegenerate as $\hat\Lu$ is 
semisimple. Also, it turns out that 
$\inn{\cdot}{\cdot}_{\mathfrak u}$ is the restriction of 
$\beta$ to $\Lu$. 

Now, since $\beta|_{W\times W}$ and 
$\inn{\cdot}{\cdot}|_{W\times W}$ are both positive-definite
real-valued symmetric bilinear 
forms which are $\Lu$-invariant, and $\tau_u$ is irreducible, there
exists $\lambda>0$ such that $\beta(x,y)=\lambda\inn xy$ 
for $x$, $y\in W$. By $\C$-bilinearity, 
$\beta^c(x,y)=\lambda\inn xy$ for $x$, $y\in V$, where $\beta^c$ is the
Killing form of $\hat{\mathfrak u}^c=\Lg+V$. 
It suffices to prove that 
$\hat{\mathfrak g}_{\mathbf R}:=\Lgr+\Vr$ is a real subalgebra of $\hat{\mathfrak u}^c$. It is clear that 
$[\Lgr,\Lgr]\subset\Lgr$ and $[\Lgr,\Vr]\subset\Vr$. 
We claim that also $[\Vr,\Vr]\subset\Lgr$.
In fact,
\begin{equation}\label{1}
 \beta^c(\ts x,\ts y)=\lambda\inn{\ts x}{\ts y}=
\lambda\overline{\inn xy}=\overline{\lambda\inn xy}
=\overline{\beta^c(x,y)} 
\end{equation}
for $x$, $y\in V$. If $Z_1$, $Z_2\in\Lg$, then also
\begin{eqnarray}\label{2} \nonumber
\beta^c(\sigma Z_1,\sigma Z_2)&=&
\mathrm{trace}_{\hat{\mathfrak u}^c}(\ad{\sigma Z_1}\ad{\sigma Z_2})\\ \nonumber
&=&\mathrm{trace}_{\mathfrak g}(\ad{\sigma Z_1}\ad{\sigma Z_2})
+\mathrm{trace}_{V}(\ad{\sigma Z_1}\ad{\sigma Z_2})\\
&=&\mathrm{trace}_{\mathfrak g}(\sigma\ad{Z_1}\ad{Z_2}\sigma)
+\mathrm{trace}_{V}(\ts\ad{Z_1}\ad{Z_2}\ts)\\ \nonumber
&=&\overline{\mathrm{trace}_{\mathfrak g}(\ad{Z_1}\ad{Z_2}})
+ \overline{\mathrm{trace}_{V}(\ad{Z_1}\ad{Z_2}})\\ \nonumber
&=&\overline{\mathrm{trace}_{\hat{\mathfrak u}^c}(\ad{Z_1}\ad{Z_2}})\\ \nonumber
&=&\overline{\beta^c(Z_1,Z_2)}, 
\end{eqnarray}
where we used in the third equality that 
\[ \ad{\sigma Z_1}x=\sigma Z_1\cdot x=\ts(Z_1\cdot\ts x)
=\ts\ad{Z_1}\ts x \]
for $x\in V$. Therefore 
\begin{align*}
\beta^c(Z,\sigma[x,y])&=\overline{\beta^c(\sigma Z,[x,y])} 
&&\text{by~(\ref{2})}\\
&=\overline{\beta^c(\sigma Z\cdot x,y)}\\
&=\beta^c(\ts(\sigma Z\cdot x),\ts y) &&\text{by~(\ref{1})}\\
&=\beta^c(Z\cdot\ts x,\ts y)\\
&=\beta^c(Z,[\ts x,\ts y])
\end{align*}
for all $Z\in\Lg$ and $x$, $y\in V$. Hence 
$\sigma[x,y]=[\ts x,\ts y]$, proving the claim.
Of course, $\rho$ is now the isotropy representation 
of the pseudo-Riemannian symmetric space
$\hat{G}_{\mathbf R}/\Gr$, where
$\hat{G}_{\mathbf R}:=\mathrm{Int}(\hat{\mathfrak g}_{\mathbf R})$ 
and $\Gr$ is the 
connected subgroup associated to $\Lgr$. \EPf

\medskip

\textit{Proof of Theorem~\ref{classif}.}
In view of Lemmas~\ref{real-polar-reducible} and~\ref{decomp}, 
it is enough to consider 
the following two cases:
\begin{enumerate}
\item $\rho$ is irreducible.
\item $\rho$ decomposes as $\rho_0\oplus\rho_0^*$, where 
$\rho_0:\Gr\to GL(V_0)$ is irreducible and non-orthogonalizable,
$\Vr=V_0\oplus V_0^*$, 
and the inner product on $\Vr$ is given by~(\ref{inner}). 
\end{enumerate}

\textit{(a.1)} Suppose first that $\rho$ is absolutely irreducible. 
Then $\tau_u$ is an absolutely irreducible polar representation 
of a compact Lie group, so it is either a Riemannian
s-representation and then the result 
follows from Lemma~\ref{irred-s-repr}, or it is listed
in~\cite{EH1}. In the latter case, it must be 
$(SO(3)\times Spin(7),\R^3\otimes\R^8)$, 
where $\R^8$ denotes the spin representation; 
according to~\cite[Table~5, p.79]{onishchik}, 
$\Gr$ equals $SO_0(1,2)\times Spin(7)$ 
(resp.~$SO(3)\times Spin_0(3,4)$, $SO_0(1,2)\times Spin_0(3,4)$;
here the subscript denotes the identity component),
and $\rho:\Gr\to GL(\R^3\otimes\R^8)$ is the tensor
product of the standard representation and the spin 
representation. 
Since $Spin(7)\subset SO(8)$
and $Spin_0(3,4)\subset SO_0(4,4)$~\cite[Thm.~14.2]{harvey},
and $\rho$ extends to a pseudo-Riemannian s-representation~$\rho'$ of 
$SO_0(1,2)\times SO(8)$ (resp.~$SO(3)\times SO_0(4,4)$,
$SO_0(1,2)\times SO_0(4,4)$) on
$\R^3\otimes\R^8$, it follows from Lemma~\ref{same-closed-orbits} that 
$\rho$ has the same closed orbits as $\rho'$, so this case is checked. 

\textit{(a.2)} Suppose now that $\rho$ is irreducible but not absolutely 
irreducible. Then $\Vr$ admits an invariant complex structure. 

\ \textit{(a.2.1)} If $\tau_u$ is irreducible, then $W$ admits 
an $U$-invariant complex structure, and by Lemma~\ref{irred-s-repr}
we have only to consider the cases in which it is not an 
s-representation. According to~\cite{EH1}, those are
\begin{equation}\label{polar-hermitian}
 \begin{array}{ll}
(SO(2)\times G_2,\R^2\otimes\R^7)& \\
(SO(2)\times Spin(7),\R^2\otimes\R^8)& \\
(SU(p)\times SU(q),(\C^p\otimes\C^q)^r)&(p\neq q) \\
(SU(n),(\Lambda^2\C^n)^r)&\mbox{($n$ odd)}\\
(Spin(10),\C^{16}) 
\end{array} 
\end{equation}
We do only the first and third cases, the others being similar in 
spirit. In the first case, 
$\Gr$ must be $SO(2)\times G_2^*$, where $G_2^*$ is the 
automorphism group of the split octonions 
and $\rho$ is the real tensor product of the standard 
representation of $SO(2)$ and the $7$-dimensional representation 
of $G_2^*$ since $\Vr$ admits an invariant complex structure.
Now $G_2^*\subset SO_0(3,4)$ and there exists an obvious 
pseudo-Riemannian s-representation $\rho'$ of $SO(2)\times SO_0(3,4)$
on $\R^2\otimes\R^7$.
It follows from Lemma~\ref{same-closed-orbits} 
that $\rho$ and $\rho'$ have the same closed 
orbits and we are done with this case. 
In the third case, viewing $\rho$ as a complex representation,
its conjugate representation $\bar\rho$ with respect to $\Gr$
must be equivalent to
$\rho^*$ because $\rho\oplus\bar\rho=(\tau_u)^c$. 
The only possibility is that $\rho$ equals
$(SU(r,p-r)\times SU(s,q-s),(\C^p\otimes\C^q)^r)$,
which has the same closed orbits as the $s$-representation of
the pseudo-Riemannian symmetric space
\[ SU(r+s,p+q-r-s)/S(U(r,p-r)\times U(s,q-s)). \]

\ \textit{(a.2.2)} If $\tau_u$ is not irreducible, then there 
exists an $U$-irreducible decomposition $W=W_1\oplus W_2$, where 
$(\tau_u)_i:U\to GL(W_i)$ is absolutely irreducible. Now 
$V=W_1^c\oplus W_2^c$ is a $G$-irreducible decomposition, where 
$W_1^c$ and $W_2^c$ are inequivalent by polarity 
(Lemma~\ref{complex-polar-reducible}(c)) and 
$W_2^c$ must be the conjugate representation to $W_1^c$ with respect 
to $\Gr$. Denote $\tau_i=(\tau_u)_i^c:G\to GL(W_i^c)$.
It follows that 
\begin{equation}\label{sigma}
\ts W_1^c = W_2^c\quad\mbox{and}\quad \tau_2(g)=\tau_1(\sigma(g))=\ts\tau_1(g)\ts
\end{equation}
for $g\in G$. Since 
$\sigma$ commutes with $\theta$, we can view $\sigma$ as an automorphism
of $U$. Suppose first that $\tau_u$ is splitting,
that is $U=U_1\times U_2$ and $\tau_u$ is the outer direct product
of $(\tau_u)_1|_{U_1}$ and $(\tau_u)_2|_{U_2}$. On the level of Lie algebras,
(\ref{sigma}) implies  
that $\Lu_1=\ker(\tau_u)_2=\sigma(\ker(\tau_u)_1)=\sigma(\Lu_2)$. 
Now we can assume that $\Lu_1=\Lu_2$, 
$(\tau_u)_1=(\tau_u)_2$,
$\Lg=\Lg_1\oplus\bar{\Lg}_1$, 
where $\Lg_1=\Lu_1^c$ and $\bar{\Lg}_1$ is 
the conjugate Lie algebra of $\Lg_1$, 
and $\sigma:\Lg\to\Lg$
is given by $\sigma(Z',\bar Z'')=(Z'',\bar Z')$. 
Moreover, $V=W_1^c\oplus \overline W_1^c$ and $\tilde\sigma\colon V\to V$ is given by $\tilde\sigma(w',\bar w'')=(w'',\bar w')$.
Hence 
$\Lgr=\{(Z',\bar Z'')\in\Lg_1\oplus\bar\Lg_1:Z'=Z''\}$,
$\Vr=\{(w',\bar w'')\in W_1^c\oplus \overline{W_1^c}:w'=w''\}$,
and $\rho$ is just the realification of the 
complexification of the real polar absolutely irreducible
representation $(\tau_u)_1|_{U_1}:U_1\to GL(W_1)$. 
If $(\tau_u)_1|_{U_1}$ is an s-representation, this means
that $\rho$ is the s-representation of a complex symmetric space
viewed as a real representation. The only other possibility is
that $(\tau_u)_1|_{U_1}$ equals 
$(SO(3)\times Spin(7),\R^3\otimes\R^8)$. In this case,
$\rho:SO(3,\C)^r\times Spin(7,\C)^r\to
GL((\C^3\otimes\C^8)^r)$ has the same closed orbits 
as $SO(3,\C)^r\times SO(8,\C)^r\to
GL((\C^3\otimes\C^8)^r)$, so we are done. 
Suppose now that $\tau_u$ is not splitting.
Then $U=U_1\times U_0\times U_2$, where $U_2$ (resp.~$U_1$)
coincides with $\ker(\tau_u)_1$ (resp.~$\ker(\tau_u)_2$) up to 
some discrete part. Since $\sigma(\ker(\tau_u)_1)=\ker(\tau_u)_2$,
the automorphism
$\sigma:U\to U$ must restrict to isomorphisms $U_1\to U_2$ and
$U_0\to U_0$. It follows that $U_0$ is essential for $(\tau_u)_1$
if and only if it is essential for $(\tau_u)_2$. Therefore
$\tau_u$ is not almost splitting in the sense 
of~\cite[p.58]{GTh1}; we use the classification 
given there: due to the facts that the $(\tau_u)_i$ admits no invariant
complex structure and $\dim W_1=\dim W_2$, we need only to 
consider the case in which 
$U_0=Spin(8)$, $U_1=U_2=\{1\}$, and $W_1$, $W_2$
are two $8$-dimensional inequivalent representations of $Spin(8)$. 
Referring to~\cite[Table~5, p.80]{onishchik}, $\Gr$ 
must be either
$Spin_0(3,5)$ or $Spin_0(1,7)$, and $\rho$ 
must be the realification
of an $8$-dimensional complex representation of $\Gr$ which is not of
real type (indeed, in each case there exist two such representations and they 
are conjugate to one another). 
 Since $\tau$ is $(SO(8,\C),\C^8_+\oplus
\C^8_-)$ (where $\C^8_{\pm}$ denote the half-spin representations)
with compact real form $(\Spin8,\R^8_+\oplus\R^8_-)$ having the 
same orbits as $(SO(8)\times SO(8),\R^8\oplus\R^8)$, 
it follows that $\rho$ has the same closed orbits as
the standard action of $(SO(8,\C))^r$
on $(\C^8)^r$. Now the latter is a pseudo-Riemannian
s-representation, so this case is also dealt with. 

\textit{(b.1)} Consider now the case in which $\rho=\rho_0\oplus\rho_0^*$,
where $\rho_0$ is absolutely irreducible and non-orthogonalizable. 
Then $\rho^c=\rho_0^c\oplus(\rho_0^*)^c$ is polar and
$\rho_0^c$ is irreducible. By polarity, $\rho_0^c$ and 
$(\rho_0^*)^c=(\rho_0^c)^*$ are inequivalent, so $\rho_0^c$ is not
self-dual implying that it is not of real type with respect to
$U$. Recall that $(\rho_0^c)^*$ is the conjugate representation
$\hat{\rho_0^c}$ with respect to $U$. 
It follows that $\tth(v',\hat v'')=(v'',\hat v')$ for 
$(v',\hat v'')\in V_0^c\oplus\hat{V}_0^c$, the space
$W=\{(v',\hat v'')\in V_0^c\oplus\hat{V}_0^c:v'=v''\}$, 
and $\tau_u:U\to GL(W)$ is irreducible, not absolutely irreducible,
and equivalent to $(\rho_0^c)^r|_U:U\to GL((V_0^c)^r)$. 
In other words, $W$ admits a $U$-invariant complex structure 
$J$ and 
$\rho_0$ is just a real form of the holomorphic extension 
of $\tau_u:U\to GL(W,J)$ to a representation of $G$ on $(W,J)$. 
By Lemma~\ref{irred-s-repr}, it suffices to consider the case
in which $\tau_u$ is not an s-representation, namely, given
in~(\ref{polar-hermitian}). We do only the case 
$(SU(n),(\Lambda^2\C^n)^r)$ for $n$ odd, the others
being similar in spirit. Since $n$ is odd,
\cite[Table~5]{onishchik} gives that $\Gr=SL(n,\R)$ and 
$\rho_0$ is the representation on $\Lambda^2\R^n$. 
Now $\rho=\rho_0\oplus\rho_0^*$ has the same closed orbits as 
$(GL^+(n,\R),\Lambda^2\R^n\oplus(\Lambda^2\R^n)^*)$, which 
turns out to be the s-representation of the pseudo-Riemannian 
symmetric space $SO(n,n)/GLâº^+(n,\R)$~\cite[Tableau~II]{berger}. 

\textit{(b.2)} Finally, suppose that 
$\rho=\rho_0\oplus\rho_0^*$,
where $\rho_0$ is irreducible, not 
absolutely irreducible and non-orthogonalizable. 
Then $\rho_0$, $\rho_0^*$ can be viewed as complex
representations, and 
$\rho^c=\rho_0\oplus\bar\rho_0\oplus\hat\rho_0\oplus\hat{
\bar{\rho_0}}$
is an irreducible decomposition with pairwise inequivalent summands,
where $\bar\rho_0$ (resp.~$\hat\rho_0=\rho_0^*$) is the 
conjugate representation to $\rho_0$ with respect to
$\Gr$ (resp.~$U$). We must have 
$\tau_u=(\tau_u)_1\oplus(\tau_u)_2:U\to GL(W_1\oplus W_2)$,
where $(\tau_u)_i$ is polar irreducible, 
not absolutely irreducible. Moreover, 
$\tau_1=\rho_0\oplus\hat\rho_0$ and
$\tau_2=\bar\rho_0\oplus\hat{\bar\rho}_0$, where we have set
$\tau_i=(\tau_u)_i^c$. 

\ \textit{(b.2.1)} Suppose $\tau_u$ is splitting. Then 
$U=U_1\times U_2$ and $\tau_u$ is the outer direct product
of $(\tau_u)_1|_{U_1}$ and $(\tau_u)_2|_{U_2}$, where 
each $(\tau_u)_i|_{U_i}$ is irreducible and not absolutely irreducible.
The automorphism $\sigma:U\to U$ must take $U_1$ to $U_2$, 
so we can assume $U_1=U_2$ and $(\tau_u)_1=(\tau_u)_2$. 
Write $G=G_1\times G_2$ where $\Lg_i=\Lu_i^c$.
Then $\rho$ is equivalent to the realification of 
$\tau_1|_{G_1}:G_1\to GL(V_0^c\oplus V_0^{c*})$, and 
$\tau_1|_{G_1}$ is the complexification of 
a polar irreducible,
not absolutely irreducible representation
$(\tau_u)_1|_{U_1}:U_1\to GL(W_1)$. We have only to consider
the case in which it is not an s-representation, namely,
given in~(\ref{polar-hermitian}). We do only the case 
$(Spin(10),(\C^{16})^r)$, the others being similar in spirit. 
Here $\tau_1$ is $(Spin(10,\C),\C^{16}\oplus\C^{16*})$ and 
$\rho$ is $(Spin(10,\C)^r,(\C^{16})^r\oplus(\C^{16*})^r)$, which 
turns out to have the same closed orbits as the pseudo-Riemannian
s-representation given by the realification
of $(\C^{\times}\times Spin(10,\C),\C^{16}\oplus\C^{16*})$.

\ \textit{(b.2.2)} Suppose $\tau_u$ is not splitting. 
Then it is not almost splitting by the same argument as
in case~(a.2.2). Owing to the fact that $(\tau_u)_i$ admits 
an invariant complex structure for $i=1$, $2$, we see 
from~\cite[p.59]{GTh1} that this case is not possible. 

\section{Isoparametric submanifolds}\label{sec:isoparam}

Let $\Vr$ be a finite-dimensional real vector space 
equipped with a non-degenerate symmetric bilinear form
$\inn{\cdot}{\cdot}$. 
A submanifold $M$ of $\Vr$ is called a \emph{pseudo-Riemannian
submanifold} if the restrictions of $\inn{\cdot}{\cdot}$ to the 
tangent spaces of $M$ are always nondegenerate.
If $M$ is a pseudo-Riemannian
submanifold, the canonical flat connection $D$ in $\Vr$ induces the 
Levi-Civit\`a connection $\nabla$ in $M$, the second fundamental form
$B$ of $M$, and the connection $\nabla^\perp$ in the normal bundle 
$\nu M$ of $M$
in the usual way. Namely,
\[ D_X Y = \nabla_X Y + B(X,Y), \]
and 
\[ D_X \xi = - A_\xi X + \nabla^\perp_X \xi, \]
where $X$ and $Y$ are sections of $TM$ and
$\xi$ is a section of $\nu M$, and the Weingarten operator
$A_\xi:TM\to TM$ is defined by
\[ \inn{A_\xi X}{Y} = \inn{B(X,Y)}{\xi}. \]
For each $p\in M$, the map
$A_\xi|_p:T_pM\to T_pM$ is a symmetric 
endomorphism with respect to the induced inner 
product in $T_pM$. Note that 
in the case in which this induced inner product is
positive-definite,
the Weingarten operator is automatically 
diagonalizable over $\R$, whereas
in the general case 
 it may happen that 
$A_\xi|_p$ is not diagonalizable, not even over $\C$.

A properly embedded
pseudo-Riemannian submanifold $M$ of $\Vr$ will be
herein called \emph{isoparametric} 
if the following two conditions are satisfied:
\begin{enumerate}
\item the normal connection is flat; 
\item the Weingarten operator
along a locally defined parallel normal vector field is
diagonalizable over $\C$ with constant eigenvalues.
\end{enumerate}
Isoparametric submanifolds of Euclidean spaces
are very important in submanifold geometry and share 
a very rich history and an extensive literature,
 see~\cite{Te,Th5,BCO} and the references therein. 
On the other hand, isoparametric submanifolds of indefinite space 
forms are not as common, but have already been considered 
before with different definitions, see e.g.~\cite{Hahn,Magid}. 
Herein we consider a stronger definition which in our opinion 
seems more natural in view of Theorems~\ref{classif} and~\ref{isoparam}. 

In this section, we will consider homogeneous isoparametric
submanifolds. We start with the following lemma. 

\begin{lem}\label{orthogonality-cartan-subspace-orbits}
Let $\tau:G\to O(V,\inn{\cdot}{\cdot})$ be a complex polar orthogonal
representation of a complex reductive algebraic group. 
\begin{enumerate}
\item For $v\in V$, we have $c_v\subset(\Lg\cdot v)^\perp$, and 
the equality holds if and only if $v$ is regular. 
\item If $c\subset V$ is a Cartan subspace, then $\inn{c}{\Lg\cdot c}=0$. 
In particular,
the restrictions of $\inn{\cdot}{\cdot}$ to $c$ and $\Lg\cdot v$ for regular $v$ 
are nondegenerate. 
\end{enumerate}
\end{lem}

\Pf (a) If $x\in c_v$, then $\Lg\cdot x\subset\Lg\cdot v$, so  
\[ \inn{x}{\Lg\cdot v}=\inn{\Lg\cdot x}v=\inn{\Lg\cdot v}v=0 \]
by $G$-invariance of $\inn{\cdot}{\cdot}$, proving the inclusion. If $v$ is 
regular, 
\begin{align*}
\dim c_v &= \dim V\cat G \\
         &= \dim V - \max_{u\in V}\dim Gu &&\text{($\tau$ is stable)} \\
         &=\dim V-\dim \Lg\cdot v,
\end{align*}
and this shows that $c_v$ is the orthocomplement of $\Lg\cdot v$ in $V$. 

(b) Follows from~(a). \EPf

\medskip

Before stating the next theorem, a couple of remarks are in order.
Let $\Gr$ be a connected real form of a connected complex reductive 
algebraic group $G$,
let $\rho:\Gr\to GL(\Vr)$ be an arbitrary real representation, and
let $\tau:G\to GL(V) $ be the complexification of $\rho$. 
If $v\in V$ is semisimple, then the isotropy subgroup
$G_v$ is reductive; hence, there exists 
a $G_v$-invariant subspace $N_v\subset V$ such that 
$V=\Lg\cdot v\oplus N_v$. The restriction of $\tau$ 
to $G_v\to GL(N_v)$ is called the \emph{slice representation} at $v$. 
If $v\in\Vr$, then 
$G_v$, $N_v$ and the slice representation are defined
over $\R$. 
There exists a Zariski-open and dense subset $V_\mathrm{pr}$
of $V$ such that all isotropy subgroups $G_v$ for semisimple 
$v\in V_\mathrm{pr}$ are in one conjugacy class~\cite[Cor.~5.6]{S}.
A semisimple point $v\in V_\mathrm{pr}$ is called 
\emph{principal}. Every principal point is regular. 
We have that $V_\mathrm{pr}\cap\Vr$ is dense in $\Vr$~\cite[13.4]{Br}.
If $v\in\Vr$, then 
$\Gr v$ is closed if and only if $Gv$ is closed~\cite{Bi},
and $\dim_{\mathbf R}\Gr v = \dim Gv$; it follows
that 
$\max_{v\in V_{\mathbf R}}\dim_{\mathbf R}\Gr v=\max_{v\in V}\dim Gv$. 
Suppose now that $\rho$ is orthogonalizable; then so is $\tau$, hence 
$\tau$ is stable; in this case, $V_\mathrm{pr}$ consists 
of semisimple elements only, and it follows from this discussion
that $V_\mathrm{pr}\cap\Vr$ is an open and dense subset of $\Vr$ 
consisting of closed $\Gr$-orbits. Suppose now, in addition, that 
$\tau$ is polar. Since the slice representations of $\tau$
are the complexifications of the slice representations
of the real form $\tau_u:U\to GL(W)$~\cite[Cor.~5.9]{S}, it follows
from~\cite[Cor.~5.4.3]{BCO} that $V_\mathrm{pr}$ is precisely the 
set of regular points of $\tau$. 

\begin{thm}\label{isoparam}
Let $\rho:\Gr\to O(\Vr,\inn{\cdot}{\cdot})$ be an orthogonal 
representation.
If $\rho$ is polar then every orbit of~$\rho$
through a regular element $v\in\Vr$ is isoparametric. 
Conversely, if $\rho$ is irreducible and 
there exists a regular element $v\in\Vr$ 
such that $\Gr v$ is isoparametric then $\rho$ is polar. 
\end{thm}

\Pf Suppose $\rho$ is polar and $v\in\Vr$ is regular. 
Then $c_v=(\Lg\cdot v)^\perp$ is a Cartan subspace of $\tau=\rho^c$
defined over $\R$. Denote the set of real points of $c_v$
by $(c_v)_{\mathbf R}$ and let $M=\Gr v$. 
Then the normal space $\nu_vM=(c_v)_{\mathbf R}$.
Since $\Lg_v\cdot c_v=0$~\cite[Lem.~2.1(iii)]{DK}, 
any $\xi\in\nu_vM$ extends to a 
locally defined equivariant normal vector field 
$\hat\xi$ along $M$ given by $\hat\xi(gv)=g\xi$ for $g\in(\Gr)^\circ$
(the connected component of the identity). 
For $X\in\Lgr$, we have that $\nabla_{X\cdot v}^\perp\hat\xi$
is the orthogonal projection in $\nu_vM$ of
$\frac{d}{dt}\big|_{t=0}(\exp tX)\xi=X\cdot\xi\in\Lgr\cdot\xi$.
Since $\Lgr\cdot\xi\subset\Lgr\cdot v$, it follows that 
$\nabla_{X\cdot v}^\perp\hat\xi=0$. This proves that 
a locally defined equivariant normal vector field along $M$ is parallel.
By taking a basis of $\nu_vM$, we get a locally defined
parallel normal frame along $\nu_vM$, which implies that 
$\nu_vM$ is flat. It is clear that the eigenvalues
of the Weingarten operator along an equivariant 
normal vector field are constant, and 
that operator is diagonalizable
over $\C$ by Example~\ref{example} below. Hence $M$ is isoparametric. 

Conversely, suppose $\rho$ is irreducible 
and there exists a regular element $v\in\Vr$
such that $M=\Gr v$ is isoparametric.
Irreducibility of $\rho$ yields that 
$M$ is full in $\Vr$, that is, not contained in a proper
affine subspace.
We first claim that a locally 
defined parallel normal vector field $\hat\xi$ along $M$ is equivariant. 
Let $U$ be a neighborhood of $v$ in $M$ where 
$\hat\xi$ is defined, and let $\hat\xi(v)=\xi$.
Suppose that $g(t)$ is a continuous
curve in $\Gr$ satisfying $g(0)=1$ and $g(t)v\in U$.
Consider the continuous curve $\xi(t)=g(t)^{-1}\hat\xi(g(t)v)$ in 
$\nu_vM$. By the isoparametric condition and the fact that
the action of $\Gr$ is orthogonal, we have that
$A_{\xi(t)}$ and $A_\xi$ have the same complex eigenvalues.
By connectedness of the domain
interval of $g(t)$ and the facts that they are diagonalizable over 
$\C$ and commute, we get that 
$A_{\xi(t)}=A_\xi$ for all $t$. 
Fullness of $M$ implies the injectivity 
of the map $\xi\mapsto A_\xi$, so $\xi(t)=\xi$ for all $t$. 
This proves the claim. Since the locally defined 
equivariant normal vector fields are parallel with respect
to the normal connection,
\[ X\cdot\xi=D_{X\cdot v}\hat\xi=-A_{\xi}(X\cdot v)+\nabla^\perp_{X\cdot v}\hat\xi
=-A_{\xi}(X\cdot v)\in\Lgr\cdot v, \]
where $\xi\in\nu_vM$ and $X\in\Lgr$. This proves that 
$\nu_vM\subset(c_v)_{\mathbf R}$. Since 
\[ \dim_{\mathbf R}\nu_vM=\dim_{\mathbf R}\Vr-\dim_{\mathbf R}M=\dim V-\dim Gv=\dim V\cat G, \]
we get that $\dim c_v=\dim V\cat G$ and hence $\tau=\rho^c$
(resp.~$\rho$) is polar. \EPf

\begin{eg}\label{example}
\em
Let $\tau:G\to O(V,\inn{\cdot}{\cdot})$ be a complex polar
orthogonal representation and fix an orthogonal real form 
$\rho:\Gr\to O(\Vr,\inn{\cdot}{\cdot})$ defined by $(\sigma,\ts)$.
In this example, we compute the Weingarten operator
of an orbit $M=\Gr v$ for a regular $v\in\Vr=V^{\tilde\sigma}$. 
Let $c$ be a $\tth$- and $\ts$-stable
Cartan subspace of $\tau$ and consider the corresponding
root space decomposition 
\[ \Lg = \Lm + \sum_{\alpha\in\mathscr A} \tilde\Lg_\alpha \]
(see subsection~\ref{roots} for the notation and terminology used
in this example). 
By Proposition~\ref{intersection} below, we may assume that
$v\in c^{\tilde\sigma}$. Let $\xi$ be a vector normal to $M$ 
at $v$ in $\Vr$. Then also $\xi\in c^{\tilde\sigma}$. 

If $\alpha$ is a noncomplex root, $\tilde\Lg_\alpha$ is 
$\sigma$-stable. 
We have (the superscript ``$\top$'' denotes the tangential 
component to the orbit)
\[ A_\xi(X_\alpha\cdot v)=-(X_\alpha\cdot\xi)^{\top}, \]
where $X_\alpha\in\tilde\Lg_\alpha^{\sigma}$, and 
\[ X_\alpha\cdot v=\alpha(v)X_\alpha\cdot v_\alpha,\quad 
X_\alpha\cdot\xi=\alpha(\xi)X_\alpha\cdot v_\alpha, \]
so
\[ A_\xi(X_\alpha\cdot v_\alpha)=\lambda\, X_\alpha\cdot v_\alpha
\quad\mbox{(resp.~$A_\xi(iX_\alpha\cdot v_\alpha)=\lambda\, 
i(X_\alpha\cdot v_\alpha)$)} \]
where $\lambda=-\frac{\alpha(\xi)}{\alpha(v)}$ is a real eigenvalue 
and $X_\alpha\cdot v_\alpha$ (resp.~$i(X_\alpha\cdot v_\alpha)$)
is the associated eigenvector if $\alpha$ is real 
(resp.~imaginary).

If $\alpha$ is a complex root, $\tilde\Lg_\alpha$ is not
$\sigma$-stable and 
$(\tilde\Lg_\alpha\oplus\tilde\Lg_{|\sigma\alpha|})^{\sigma}$
is spanned by $X_\alpha+\sigma X_\alpha$ and $i(X_\alpha-\sigma X_\alpha)$ 
for $X_\alpha\in\tilde\Lg_\alpha^\theta$. The associated
subspace of $T_vM$ is spanned by
\begin{equation}\label{basis}
\alpha(v)X_\alpha\cdot v_\alpha
+\overline{\alpha(v)}\ts(X_\alpha\cdot v_\alpha),\quad
i\left(\alpha(v)X_\alpha\cdot v_\alpha
-\overline{\alpha(v)}\ts(X_\alpha\cdot v_\alpha)\right) 
\end{equation}
for $X_\alpha\in\tilde\Lg_\alpha^\theta$. 

Now 
$\lambda=-\frac{\alpha(\xi)}{\alpha(v)}$
is not real and the matrix of $A_\xi$ in the basis~(\ref{basis})
is given by 
\[ \left(\begin{array}{cc}
                     \Re\lambda & -\Im\lambda \\
                     \Im\lambda & \Re\lambda
       \end{array} \right), \]
which is of course diagonalizable over~$\C$.
\end{eg}

\section{Structural theory of 
polar representations of\\ real reductive algebraic 
groups}\label{sec:structure}

Consider a semisimple pseudo-Riemannian symmetric 
space $\hGr/\Gr$ and its complexification $\hat G/G$ as
in the first two paragraphs of section~\ref{sec:classif}. 
Let $\hat\sigma$ denote the conjugation of $\hat G$ 
over $\hat{G}_{\mathbf R}$. 
We can choose a Cartan involution $\hat\theta$ of $\hGr$ 
that commutes with $\hat\tau$ on $\hGr$. Since $\hat G$  
is simply-connected, we can extend $\hat\theta$ anti-holomorphically
to a Cartan involution of $\hat G$ which will be denoted by the same
letter. Note that $\hat\theta$ commutes with $\hat\tau$ and $\hat\sigma$ 
on $\hat G$. Set $\theta$ (resp.~$\tth$)
to be the restriction of $\hat\theta$ to $G$ (resp.~$V$), and set 
$\sigma$ (resp.~$\ts$)
to be the restriction of $\hat\sigma$ to $G$ (resp.~$V$). 
Then $\hat U=\hat G^{\hat\theta}$ (resp.~$U=G^\theta$)
is a compact real form of $\hat G$ (resp.~$G$). Write
$W=V^{\tilde\theta}$.
Now we have the combined decomposition
\begin{equation}\label{combined}
 \Lhgr = (\underbrace{\Lgr\cap\Lu}_{:=\Lkr}\stackrel{\theta}{+}
\underbrace{\Lgr\cap i\Lu}_{:=\Lpr})
\stackrel{\hat\tau}{+}(\Vr\cap W\stackrel{\tilde\theta}{+}\Vr\cap iW). 
\end{equation}
In this context, an element $v\in\Vr$ is called
semisimple if $\mathrm{ad}_v$ is a semisimple endomorphism 
of $\Lhg$, and a Cartan subspace of $\Lhgr$ is a maximal
Abelian subspace of $\Vr$ consisting of semisimple elements.  
It is known that the $\mathrm{Ad}_{G_{\mathbf R}}$-orbit
of $v\in \Vr$ is closed if and only if $v$ is 
semisimple~\cite[Cor.~10.3]{BH};
every semisimple element of $\Vr$ belongs to some Cartan subspace;
every Cartan subspace of $\Vr$ is 
$\mathrm{Ad}_{(G_{\mathbf R})^\circ}$-conjugate
to a $\tth$-stable Cartan subspace; there exist finitely 
many $\mathrm{Ad}_{(G_{\mathbf R})^\circ}$-conjugacy classes of 
$\tth$-stable Cartan subspaces in $\Vr$; two 
such $\tth$-stable Cartan subspaces are 
$\mathrm{Ad}_{(K_{\mathbf R})^\circ}$-conjugate
if and only they are 
$\mathrm{Ad}_{(G_{\mathbf R})^\circ}$-conjugate if and only 
they are $\mathrm{Ad}_{G}$-conjugate~\cite{HHNO}.

Throughout this section, we let 
$\tau:G\to GL(V)$ be a complex polar 
representation of a connected complex reductive algebraic 
group, consider a real form 
$\rho:\Gr\to GL(\Vr)$ defined by $(\sigma,\ts)$,
where $\Gr$ is the identity component of $G^\sigma$,
and prove a collection
of results for $\rho$ similar to those stated in the previous 
paragraph for an s-representation. 
The first three results do not require
that~$\tau$ and~$\rho$ be orthogonalizable. 

\subsection{General facts about Cartan subspaces}

A \emph{Cartan subspace} of $\rho$ is a subspace 
of $V^{\tilde\sigma}$ which is the $\ts$-fixed point vector space 
of a $\ts$-stable Cartan subspace of $\tau$.

\begin{lem}\label{ts-stable-cartan-subspaces}
There exist $\ts$-stable Cartan subspaces of $\tau$.
\end{lem}

\Pf Owing to the remarks preceding Theorem~\ref{isoparam},
the set $V_{\mathrm{pr}}\cap\Vr$ is a nonempty open subset 
of $\Vr=V^{\tilde\sigma}$ consisting of regular elements of $\tau$;
it suffices to take $c_v$ where $v$ lies therein. \EPf 

\medskip

We will use the following notion in the proof
of the next proposition. The \emph{rank} of $\tau$ 
is defined to be the difference $\dim c-\dim c^{\mathfrak g}$,
where $c\subset V$ is a Cartan subspace and 
$c^{\mathfrak g}$ denotes the subspace of $G$-fixed points in~$c$. 

\begin{prop}\label{intersection}
Given a semisimple $x\in V^{\tilde\sigma}$, there exists a 
Cartan subspace of $V^{\tilde\sigma}$ which contains $x$.
\end{prop}

\Pf Note that for a regular $x\in V^{\tilde\sigma}$, one can simply 
take $c=c_x$. In the general case, we proceed 
by induction on the rank of $\tau$. Since $x$ is semisimple,
there exists a Cartan subspace $c'$ such that $x\in c'$. 
If $x\in (c')^{\mathfrak g}$, then $x$ belongs to any $\ts$-stable
Cartan subspace of~$\tau$.
Suppose now $x\not\in (c')^{\mathfrak g}$. Then the slice representation
$(G_x,N_x)$ is polar with rank stricly lower 
than $\tau$, and $c'\subset N_x$ is a Cartan subspace 
of $(G_x,N_x)$~\cite[Thm.~2.4]{DK}. Without loss of generality,
$x$ is a minimal vector with respect to some $U$-invariant 
positive-definite
Hermitian form $(\cdot,\cdot)$ which is real-valued on $V^{\tilde\sigma}$, 
and $N_x$ is the orthocomplement of $\Lg\cdot x$
with respect to $(\cdot,\cdot)$~\cite[Rmk.~1.4]{DK}. Since $x\in V^{\tilde\sigma}$, it 
follows that $G_x$ is $\sigma$-stable, $N_x$ is 
$\ts$-stable and $(G_x,N_x)$ is defined over $\R$ with respect
to $(\sigma,\ts)$. By the induction hypothesis, there exists a
$\ts$-stable Cartan subspace $c\subset N_x$ such that 
$x\in c$. Now $c$, $c'$ are two Cartan subspaces of
$(G_x,N_x)$, so there exists $g\in G_x$ such that $g\cdot c'=c$.
It follows that $c$ is a Cartan subspace of $\tau$. \EPf 

\begin{thm}\label{finiteness}
There exist only finitely many $\Gr$-conjugacy classes 
of Cartan subspaces of~$\Vr$.
\end{thm}

\Pf According to the remarks preceding Theorem~\ref{isoparam},
the set of regular points of $\tau$ is a Zariski-open
and dense subset $V_{\mathrm{pr}}$ of $V$. By a theorem of 
Whitney~\cite{W}, $V_{\mathrm{pr}}\cap\Vr$
has finitely many connected components. 

Suppose now that $c^{\tilde\sigma}$ is a Cartan subspace 
of~$\Vr$. Consider the map 
\[ \Gr\times c^{\tilde\sigma}\to\Vr,\qquad
(g,v)\mapsto g\cdot v; \]
it is easily seen to be a smooth submersion at $v$ 
if $v$ is a regular point
of $\tau$. It follows that 
$\Gr\cdot(c^{\tilde\sigma}\cap V_{\mathrm{pr}})$ is open 
in $\Vr$. But the sets 
$\Gr\cdot(c^{\tilde\sigma}\cap V_{\mathrm{pr}})$
for varying $c^{\tilde\sigma}$ obviously cover 
$V_{\mathrm{pr}}\cap\Vr$. 
Any two of them are not disjoint 
if and only if the corresponding Cartan subspaces are conjugate, in 
which case the sets coincide. The result follows. \EPf

\medskip

Consider the categorical quotient map $\pi:V\to V/\!/G$.
Since $G$, $V$, and the action of $G$ on $V$ are defined 
over $\R$, so is the variety $V/\!/G$; denote its set of real
points by $(V/\!/G)_{\mathbf R}$. By a theorem of Tarski 
and Seidenberg, $\pi(\Vr)$ is a real semialgebraic subset 
of~$(V/\!/G)_{\mathbf R}$. Recall that 
$\pi(V_\mathrm{pr}\cap\Vr)$ is an open and dense 
subset of $\pi(\Vr)$. We propose the following
conjecture (compare~\cite{Ro}). 

\begin{conj}
The map $\pi$ sets up a one-to-one correspondence 
between the $\Gr$-conjugacy classes of Cartan 
subspaces of $\Vr$ and the connected components of the 
stratum $\pi(V_\mathrm{pr}\cap\Vr)$.
\end{conj}

Henceforth we assume that $\rho$ (and hence $\tau$) is orthogonal
with respect to~$\inn{\cdot}{\cdot}$. 

\begin{thm}\label{conjugacy-cartan-subspaces}
\begin{enumerate}
\item Given a $\ts$-stable Cartan subspace $c\subset V$, there 
exists a Cartan pair $(\eta,\tilde\eta)$ commuting with $(\sigma,\ts)$
such that $c$ is $\tilde\eta$-stable.
\item Given a Cartan pair $(\theta,\tth)$ commuting with $(\sigma,\ts)$,
every $\ts$-stable Cartan subspace $c\subset V$ is 
$(G^\sigma)^\circ$-conjugate 
to a $\tth$-stable one (hence also $\ts$-stable).
\end{enumerate}
\end{thm}

\Pf We begin by showing that there exists a 
Cartan pair $(\mu,\tilde\mu)$ of $\tau$ 
such that $\tilde\mu(c)=c$. 
Indeed, suppose $(\theta',\tilde\theta')$ is any Cartan pair.
We can select $v\in V^{\tilde\theta'}$ regular. Since $c$ meets 
all the closed orbits, there exists $g\in G$ such that 
$g\cdot v\in c$. 
Define $\mu=\mathrm{Inn}_g\theta'\mathrm{Inn}_g^{-1}$ 
and~$\tilde\mu=g\tth'g^{-1}$. Then $(\mu,\tilde\mu)$ is a Cartan pair
and $\tilde\mu(g\cdot v)=g\cdot v$. Hence $c=c_{g\cdot v}$ is 
$\tilde\mu$-stable. 

The following construction of $\eta$ is standard 
(compare~\cite[\S3, Prop.~6]{onishchik}).
Set $\omega=\sigma\mu$. We can view $\omega$ as a complex
linear automorphism of $\Lg$. Consider the decomposition 
into the center and semisimple factor $\Lg=\Lz\oplus\Lg_{ss}$.
Let $\beta$ be the Killing form of $\Lg_{ss}$. 
We can extend $\beta$ to an $\ad{}$-invariant symmetric bilinear form 
on $\Lg$, denoted by the same letter, which is real-valued
on $\Lg^{\mu}$, $\Lg^{\sigma}$ and negative-definite 
on $\Lg^{\mu}$. Then one easily sees that 
$\omega$ is Hermitian
with respect to the positive-definite
Hermitian form $B_{\mu}(X,Y)=-\beta(X,\mu Y)$, where $X$, $Y\in\Lg$. 
It follows that 
$\omega^2$ is Hermitian and 
positive-definite, and hence belongs to a one-parameter
family of Hermitian and positive-definite automorphisms of $\Lg$.
Therefore there exists a 
unique Hermitian, positive-definite automorphism $\varphi$ of 
$\Lg$ such that $\varphi^4=\omega^2$. 
Since $\varphi|_{\mathfrak g_{ss}}$ belongs to a one-parameter 
group of automorphisms of 
$\Lg$, we have that $\varphi|_{\mathfrak g_{ss}}$ is inner, 
that is, equals $\Ad{h}$ for some $h\in G_{ss}$.
Set $\eta=\mathrm{Inn}_h\,\mu\,\mathrm{Inn}_h^{-1}$. 
Then $\eta$ is a Cartan involution
of $G$. Also, on the Lie algebra level,
$\mu\omega\mu=\omega^{-1}$, so 
$\mu\omega^2\mu=\omega^{-2}$ and $\mu\varphi\mu=\varphi^{-1}$.
Of course, $\omega\omega^2\omega^{-1}=\omega^2$, so
$\omega\varphi\omega^{-1}=\varphi$ and 
$\omega\varphi^2\omega^{-1}=\varphi^2$. Now we have
\[ \eta\sigma=\varphi\mu\varphi^{-1}\sigma=\varphi^2\mu\sigma=\varphi^2
\omega^{-1}=\omega^{-1}\varphi^2, \]
\[ \sigma\eta=\sigma\varphi\mu\varphi^{-1}=\sigma\mu\varphi^{-2}=\omega\varphi^{-2}, \]
so $\varphi^4=\omega^2$ implies that $\eta\sigma=\sigma\eta$
on $\Lg$, and also on $G$.

For the next step, define $\tom=\ts\tmu$. Then $\tom$ is a 
$G$-equivariant 
complex automorphism of $V$. 
Futher, $\tom$ is Hermitian with respect
to the positive-definite Hermitian form 
$B_{\tilde\mu}(x,y)=-\inn x{\tmu y}$ on $V$. It follows
that $\tom^2$ is Hermitian and positive-definite, so as above
there is a unique Hermitian and positive-definite automorphism 
$\tilde\varphi$ of $V$ such that 
$\tilde\varphi^4=\tom^2$.
Setting 
$\tilde\eta=\tilde\varphi\tmu\tilde\varphi^{-1}$, we have that
$\tilde\eta\tilde\sigma=\tilde\sigma\tilde\eta$ by a computation
similar to that in the previous paragraph. Moreover, $\tilde\eta(c)=c$,
because $\ts(c)=c$ and $\tmu(c)=c$. We also have ($x$, $y\in V$)
\begin{eqnarray*}
\inn{\tilde\eta x}{\tilde\eta y} & = & 
\inn{\tilde\varphi\tmu\tilde\varphi^{-1}x}{\tilde\varphi\tmu\tilde\varphi^{-1}y} \\
&=& \inn{\tmu\tilde\varphi^{-1}x}{\tmu\tilde\varphi^{-1}y} \\ 
&=&\overline{\inn{\tilde\varphi^{-1}x}{\tilde\varphi^{-1}y}} \\ 
&=&\overline{\inn xy}
\end{eqnarray*}
and, if $0\neq x\in V^{\tilde\eta}$,
\begin{align*}
\inn xx&= \inn{\tilde\varphi^{-1}x}{\tilde\varphi^{-1}x} \\
     &<0, &&\text{($\tilde\varphi^{-1}x\in V^{\tilde\eta}$)} 
\end{align*}
where we have used that ($x$, $y\in V$)
\begin{align*}
 \inn{\tilde\varphi x}{\tilde\varphi y}& = 
-B_{\tilde\mu}(\tilde\varphi x,\tmu\tilde\varphi y) \\
&= -B_{\tilde\mu}(\tilde\varphi x,\tilde\varphi^{-1}\tmu y) 
&&\text{($\tmu\tilde\varphi\tmu=\tilde\varphi^{-1}$)}\\
&= -B_{\tilde\mu}(x,\tmu y)  
&&\text{($\tilde\varphi$ is Hermitian)}\\
&=  \inn xy\\
&=\inn{\tilde\varphi^{-1}x}{\tilde\varphi^{-1}y}.
\end{align*}

In order to see that $(\eta,\tilde\eta)$ is a Cartan pair,
it only remains to check that 
$\tilde\eta(g\cdot v)=\eta(g)\cdot\tilde\eta(v)$ 
for $g\in G$, $v\in V$. It suffices to prove that 
$\tilde\varphi=\tau(h)$. Denote the induced representation 
$\Lg\to\mathfrak{gl}(V)$ by $d\tau$. 
%First, note that
%$B_{\mu}(\varphi X,Y)=B_{\mu}(X,\varphi Y)$ for $X$, $Y\in\Lg$
%implies that 
%\[ \beta(X,\Ad{h^{-1}}\mu Y)=\beta(\Ad h X,\mu Y)=\beta(X,\mu\Ad h Y) \]
%so $\Ad{h^{-1}}=\Ad{\mu(h)}$ implying that $\mu(h)=h^{-1}$.
%Since $\mu$ is a Cartan involution of $G$, there exists
%$Y\in i\Lg^{\mu}$ such that $\exp Y=h$. 
Since $\Ad h$ is Hermitian, positive-definite with respect to 
$B_{\mu}$, the element~$h$ 
can be taken of the form $\exp Y$, where 
$Y\in i\Lg_{ss}^\mu$. 
Then $\tau(h)=e^{d\tau(Y)}$. This implies that $\tau(h)$ 
is  Hermitian, positive-definite with respect to $B_{\tilde\mu}$. 
Since ($X\in\Lg$) 
\[ \ts d\tau(X)\ts^{-1} = d\tau(\sigma X)\quad\mbox{and}\quad
\tmu d\tau(X)\tmu^{-1} = d\tau(\mu X), \]
we also have that
\begin{equation}\label{omega}
 \tom d\tau(X)\tom^{-1} = d\tau(\omega X). 
\end{equation}
Since the irreducible summands of $V$ must be pairwise 
inequivalent by polarity, each one of them is $\tom$-invariant. 
Let $V_0$ be an irreducible summand of $V$ and suppose
that the action of $\Lz$ on $V_0$ is given by a linear
functional $\Lambda:\Lz\to\C$. Equation~(\ref{omega})
implies that $\Lambda(X)=\Lambda(\omega X)$ for $X\in\Lz$. 
Now, if $X\in\Lz$ and $v\in V_0$, we have
\begin{equation}\label{tau-h-1}
 \tau(h)^4d\tau(X)\tau(h)^{-4}v=\Lambda(X)v=\Lambda(\omega^2(X))v
=d\tau(\omega^2(X))v, 
\end{equation}
and if $X\in\Lg_{ss}$,
\begin{eqnarray}\label{tau-h-2}
 \tau(h)^4d\tau(X)\tau(h)^{-4}
&=& d\tau(\mathrm{Ad}_h^4 X) \nonumber \\ 
&=& d\tau(\varphi^4(X)),\\ 
&=& d\tau(\omega^2(X)). \nonumber
\end{eqnarray}

Equations~(\ref{omega}), 
(\ref{tau-h-1}) and~(\ref{tau-h-2}) imply
that $\tom^2$ and $\tau(h)^4$ are two intertwining
maps between the representations $d\tau$ and $d\tau\circ\omega^2$.
It follows that they are multiples of each other on each 
irreducible summand. Since both maps
are positive-definite, $\tau(h)^4=\lambda \tom^2$ for 
$\lambda\in\R$, $\lambda>0$. Since both are isometries with respect to
$\inn{\cdot}{\cdot}$, one has 
$\lambda=1$. Now~(a) is proved. For proving~(b), 
construct $(\eta,\tilde\eta)$ as in~(a) and note that 
it is conjugate to 
$(\theta,\tth)$ by an element $g'\in(G^\sigma)^\circ$ by 
Corollary~\ref{conjugate-cartan-pairs}. 
Now $c'=g'\cdot c$ is a $\tth$-stable Cartan 
subspace. \mbox{}\hfill\EPf

\medskip

In case 
a Cartan pair $(\theta,\tth)$ commuting with 
$(\sigma,\ts)$ is fixed, a $\tth$-stable 
Cartan subspace of $\rho$ will sometimes
be called \emph{standard}.

\begin{cor}\label{cor:intersection}
If $(\theta,\tth)$ is a Cartan pair 
commuting with $(\sigma,\ts)$, then every closed
$(G^{\sigma})^\circ$-orbit in $V^{\tilde\sigma}$ intersects a
standard Cartan subspace of $ V^{\tilde\sigma}$.
\end{cor}

\Pf Suppose that $(G^{\sigma})^\circ x$ is a closed orbit in
$V^{\tilde\sigma}$. By Proposition~\ref{intersection}, 
there exists a $\ts$-stable Cartan subspace $c\subset V$ 
such that $x\in c^{\tilde\sigma}$. 
By Theorem~\ref{conjugacy-cartan-subspaces}, there exists 
$g\in(G^\sigma)^\circ$ such that $g\cdot c$ is a 
$\ts$- and $\tth$-stable Cartan subspace. Of course, 
$(G^\sigma)^\circ x$ meets $g\cdot c$. \EPf

\subsection{Roots and co-roots}\label{roots}

In the rest of the paper, we assume that $\rho$ is orthogonal 
with respect to~$\inn{\cdot}{\cdot}$ and a 
Cartan pair $(\theta,\tth)$ commuting with $(\sigma,\ts)$
has been fixed according to Proposition~\ref{cartan-pairs}.
We also recall the Hermitian form $(\cdot,\cdot)$ that was introduced 
in that proposition and satisfies 
equation~(\ref{hermitian-symmetric}).   

For a given Cartan subspace $c\subset V$, the set
of singular elements  
$c_{\mathrm{sing}}\subset c$ is by definition
the complement of the set of regular elements in~$c$.
If the rank of $\tau$ is not zero, it is known that 
$c_{\mathrm{sing}}$ is a union of finitely many complex hyperplanes 
\[ c_{\mathrm{sing}} = \bigcup_{\alpha\in\mathscr A} c_\alpha, \]
where $\mathscr A$ is a finite index set~\cite[Lem.~2.11]{DK}. 
Fix a $\ts$- and $\tth$-stable Cartan subspace $c\subset V$, set 
$\Lg_\alpha$ to be the centralizer of $c_\alpha$ in $\Lg$
and $G_\alpha$
to be the corresponding connected subgroup of $G$.

\begin{lem}\label{orthogonality-root-subspaces}
We have that $\inn{\Lg_\alpha\cdot c}{\Lg_\beta\cdot c}=0$ for 
$\alpha\neq\beta$. 
\end{lem}

\Pf It follows from Lemma~\ref{orthogonality-cartan-subspace-orbits}
that $c\subset(\Lg\cdot v)^\perp$ for $v\in c$. 
Since $(\Lg\cdot v)^\perp$ is $\Lg_v$-invariant, this 
implies
$\inn{\Lg_v\cdot c}{\Lg\cdot v}=0$.
In particular, if $v\in c_\alpha\setminus\cup_{\beta\neq\alpha}c_\beta$,
then $\Lg_v=\Lg_\alpha$~\cite[p.~516]{DK}, so
$\inn{\Lg_\alpha\cdot c}{\Lg \cdot v}=0$ implying that 
$\inn{\Lg_\alpha\cdot c}{\Lg\cdot  c_\alpha}=0$. Since
$\Lg_\beta\cdot c\subset\Lg\cdot c_\alpha$ for 
$\alpha\neq\beta$~\cite[p.~517]{DK}, the desired result follows. 
\EPf 

\begin{lem}\label{c-alpha-theta-tilde-stable}
Each $c_\alpha$ meets $V^{\pm\tilde\theta}$ in a real hyperplane.
\end{lem}

\Pf It is equivalent to prove that each $c_\alpha$ is 
$\tth$-stable.	
Of course, $(\cdot,\cdot)$ is nondegenerate on $c_\alpha\times c_\alpha$ 
as $(\cdot,\cdot)$ is positive-definite. Choose $v_\alpha\in c$ to be
$(\cdot,\cdot)$-orthogonal to $c_\alpha$. We claim that the decomposition
$c=c_\alpha\oplus\C v_\alpha$ is $\inn{\cdot}{\cdot}$-orthogonal. Since
$\inn{x}{\tth y}=-(x,y)$ for $x$, $y\in V$, this will prove 
the desired result. In order to prove the claim, note that 
$c\oplus \Lg_\alpha\cdot c$ is a $G_\alpha$-invariant 
subspace~\cite[Thm.~2.12(ii)]{DK}
and $\inn{\cdot}{\cdot}$ is nondegenerate on $c\oplus \Lg_\alpha\cdot c$
by Lemmas~\ref{orthogonality-cartan-subspace-orbits}
and~\ref{orthogonality-root-subspaces}. Since $G_\alpha$ acts 
trivially on $c_\alpha$ and 
$\Lg_\alpha\cdot c=\Lg_\alpha\cdot v_\alpha$, 
it follows that 
$c_\alpha\oplus\C v_\alpha\oplus\Lg_\alpha\cdot v_\alpha$
is $\inn{\cdot}{\cdot}$-orthogonal. \EPf

\medskip

Since $\inn{\cdot}{\cdot}$ is positive-definite on 
$V^{-\tilde\theta}\times V^{-\tilde\theta}$, the vector $v_\alpha$
in the proof of Lemma~\ref{c-alpha-theta-tilde-stable} can be chosen
to satisfy
\[ v_\alpha\in V^{-\tilde\theta},\quad
\inn{v_\alpha}{c_\alpha^{-\tilde\theta}}=0\quad\mbox{and}
\quad\inn{v_\alpha}{v_\alpha}=1, \]
and then it is uniquely defined up to a sign.
We select a connected 
component of $c^{-\tilde\theta}
-\cup_{\alpha\in\mathscr{A}} c_\alpha^{-\tilde\theta}$ 
once and for all, 
and then $v_\alpha$ is uniquely defined
(but the 
sign of $v_\alpha$ will not actually matter for our purposes). 
The vector $v_\alpha$ is called a \emph{(unnormalized) co-root}.
The associated \emph{root} is the linear functional 
$\alpha:c\to\C$ obtained by setting
\[ \alpha(v) = \inn{v}{v_\alpha}\in\R \]
for $v\in c^{-\tilde\theta}$ and then considering its 
complex-linear extension to $c$. A root is called:
\emph{real} (resp.~\emph{imaginary})
if $\alpha$ is real-valued (resp.~purely imaginary-valued) 
on $c^{\tilde\sigma}$, and it is called \emph{complex} otherwise. 
It follows from the $\inn{\cdot}{\cdot}$-orthogonality of
the decomposition 
$c^{-\tilde\theta}=\cs\cap c^{-\tilde\theta}\oplus 
c^{-\tilde\sigma}\cap c^{-\tilde\theta}$
that $\alpha$ is real (resp.~imaginary)
if and only if it vanishes on~$c^{-\tilde\sigma}\cap c^{-\tilde\theta}$ (resp.~$\csmt$), in which
case $v_\alpha$ belongs to $c^{\tilde\sigma}\cap c^{-\tilde\theta}$
(resp.~$c^{-\tilde\sigma}\cap c^{-\tilde\theta}$). 
It follows that $\alpha$ in noncomplex if and only if 
$c_\alpha$ is $\ts$-invariant if and only if 
it is $\tom$-invariant, where $\tom=\ts\tth=\tth\ts$. 
Recall 
that $\tth$ gets replaced by its opposite by changing the 
sign of $\inn{\cdot}{\cdot}$, so the choice of some signs above 
does not have intrinsic meaning, as compared to the case 
of an s-representation in which the sign of $\inn{\cdot}{\cdot}$ 
is fixed by the Killing form of $\hat{\Lg}_{\mathbf R}$ 
(see~(\ref{combined})). 

Let $\Lm$ be the centralizer of $c$ in $\Lg$. Then $\Lm$ is 
$\sigma$-, $\theta$-stable. 
Since $\Lm$ is a reductive subalgebra of $\Lg_{\alpha}$, there exists a 
$\theta$- and $\mathrm{ad}_{\mathfrak m}$-stable splitting
\[ \Lg_\alpha = \Lm\oplus \tilde\Lg_\alpha, \]
where $\tilde\Lg_\alpha$ is a subspace, which is 
called a \emph{root space}. 
Now assume $\alpha$ is noncomplex. Then 
$\tilde\Lg_\alpha$ can be taken $\omega$-stable, so that 
$\tilde\Lg_\alpha=\tilde\Lg_\alpha^{\omega}\oplus
\tilde\Lg_\alpha^{-\omega}$. An imaginary
root $\alpha\in \mathscr A$ 
is called \emph{noncompact imaginary} 
if $\tilde\Lg_\alpha^{-\omega}\neq0$
and \emph{compact imaginary} otherwise. 
A real root $\alpha\in \mathscr A$ 
is called \emph{compact real} 
if $\tilde\Lg_\alpha^{\omega}\neq0$
and \emph{noncompact real} otherwise. 
Finally, define 
\[ \ts\alpha(v)=\overline{\alpha(\ts v)}, \]
where $v\in c$. Since $\ts$ takes singular orbits to singular 
orbits and maps hyperplanes of $c$ to 
hyperplanes of $c$, this defines an action on 
$\mathscr A\cup(-\mathscr A)$. 
Also, $\ts\alpha=\alpha$ (resp.~$\ts\alpha=-\alpha$)
if and only if $\alpha$ is real (resp.~imaginary). 
We can choose the root spaces so that 
$\sigma\tilde\Lg_\alpha=\tilde\Lg_{|\tilde\sigma\alpha|}$
for all $\alpha\in\mathscr A$, where 
$|\cdot|:\mathscr A\cup(-\mathscr A)\to\mathscr A$ has 
its obvious meaning. 

\subsection{Cayley transforms}\label{cayley}

By Corollary~\ref{cor:intersection}, 
every closed $\Gr$-orbit in $\Vr$ meets some 
standard Cartan subspace of $\Vr$. 
We want to study standard Cartan subspaces of $\Vr$, so 
consider a $\ts$- and $\tth$-stable Cartan subspace $c\subset V$.
Note that 
\[ \cs = \cst \oplus \csmt, \]
and $\dim_{\mathbf R} \cst$ (resp.~$\dim_{\mathbf R} \csmt$)
is an invariant of the $\Gr$-conjugacy class
of $\cs$, called the \emph{compact dimension} 
(resp.~\emph{noncompact dimension}) of $\cs$.
We call a standard Cartan subspace 
$\cs$ \emph{maximally compact} 
(resp.~\emph{maximally noncompact}) if 
its compact dimension (resp.~\emph{noncompact dimension})
is as large as possible. Note that the compact and 
noncompact dimensions of $\cs$ are interchanged 
if we replace $\inn{\cdot}{\cdot}$ and $\tth$ by their 
opposites. 
A maximally compact or maximally noncompact standard 
Cartan subspace will also be called
\emph{extremal}.
Cayley transforms are used to pass
from one $\Gr$-conjugacy class of Cartan subspaces to another one, 
namely, 
to increase or decrease its compact dimension by one
(and correspondingly
decrease or increase its noncompact dimension by one).  
In general, an element $g\in G$ maps a $\ts$- and $\tth$-stable 
Cartan subspace $c$ of $V$ to another $\ts$- and $\tth$-stable 
Cartan subspace if and only if $\sigma(g)g^{-1}$ and 
$\theta(g)g^{-1}$
belong to the normalizer $N_G(c)$ of $c$ in $G$, 
as is easily seen. Recall that the 
\emph{Weyl group} of $c$ is the finite group~\cite[p.~513]{DK}
\[ W(c) = N_G(c)/Z_G(c), \]
where $Z_G(c)$ denotes the centralizer of $c$ in $G$.
We will construct a special
kind of Cayley transform. We first consider the case of a 
rank one polar 
orthogonal irreducible
representation $\tau:G\to O(V,\inn{\cdot}{\cdot})$. 
Fix a 
standard Cartan subspace $c$ which is extremal, say
maximally compact. Here $\cs=\cst$ and $\mathscr A=\{\alpha\}$. 
Assume 
that $\alpha$ is an imaginary root. We will show how one can pass
from $\cs$ to a Cartan subspace $\hat c$ in another 
$\Gr$-conjugacy class which in this case, by dimensional reasons, 
must be maximally noncompact, namely, $\hat c^{\tilde\sigma}=
\hat c^{\tilde\sigma}\cap\hat c^{-\tilde\theta}$.
Since the rank is one, 
$\tau_u:U\to O(W,\inn{\cdot}{\cdot})$ is a co-homogeneity one action of a 
compact Lie group. Let $v=iv_\alpha\in\cst$. 
Then $\inn{v}{v}=-1$ and 
$U(v)$ is a round sphere $S^{n-1}\approx U/U_v$ in~$W$. 
Introduce the following 
notation: $\Lg=\Lk+\Lp$ is the decomposition into 
$\pm1$-eigenspaces of 
$\omega$, $\Lk_{\mathbf R}=\Lk^\sigma=\Lk^\theta$, 
$\Lp_{\mathbf R}=\Lp^\sigma$, 
and $K_{\mathbf R}=U^\sigma=U\cap\Gr$; 
note that $K_{\mathbf R}$ is a maximal compact subgroup 
of $\Gr$ and hence it is connected since $\Gr$ is so. 

\begin{claim}
We have that $\alpha$ is compact imaginary if and only if $K_{\mathbf R}\subset U$ 
is transitive on $S^{n-1}$.
\end{claim}

In fact, here
we have $\Lg_\alpha=\Lm\oplus\tilde\Lg_\alpha$ where $\Lg_\alpha=\Lg$,
$\Lm=\Lg_v$, and $\Lg=\Lg_v\oplus\tilde\Lg_\alpha$ is 
$\theta$-stable. 
Taking $\theta$-fixed points, we get 
$\Lu=\Lu_v\oplus\tilde\Lg_\alpha^\theta$.
Now $K_{\mathbf R}$ is transitive on $S^{n-1}$ if and only if 
$\Lu_v+\Lk_{\mathbf R}=\Lu$ if and only if 
$\Lk_{\mathbf R}\supset\tilde\Lg_\alpha^\theta$ if and only if 
$\Lk\supset\tilde\Lg_\alpha$ if and only if 
$\tilde\Lg_\alpha^{-\omega}=\{0\}$. 

\begin{claim}
If $K_{\mathbf R}$ is not transitive on $S^{n-1}$, then we can take 
$g\in G$ such that $\theta(g)=g$, $\sigma(g)=g^{-1}$ and 
$g^2=-\mathrm{id}\in W(c)$. 
\end{claim}

Indeed, the assumption is equivalent to 
$\tilde\Lg_\alpha^\theta\cap\tilde\Lg_\alpha^{-\sigma}\neq\{0\}$;
take a nonzero $X$ therein. We can choose $X$ so that
$\gamma(t)=\exp tX\cdot v$ is a 
unit speed geodesic
of $S^{n-1}$ connecting $\gamma(0)=v$ to $\gamma(\pi)=-v$. Set 
$g=\exp\frac{\pi}2X\in U$. Clearly, $\theta(g)=g$. Also, 
$\sigma(g)=g^{-1}$, and $g^2\cdot v=\exp\pi X\cdot v=-v$, so 
$g^2=-\mathrm{id}$
on $\C v=c$. 

\begin{claim}
If $g$ is as in the previous claim and $\hat c = g\cdot c$, then 
${\hat c}^{\tilde\sigma}$ 
is a maximally noncompact Cartan subspace of $V^{\tilde\sigma}$. 
\end{claim}

In fact, $\theta(g)g^{-1}=\mathrm{id}$ and $\sigma(g)g^{-1}=g^{-2}=-\mathrm{id}$
both belong to $W(c)$, so $\hat c$ is $\ts$- and $\tth$-stable. Also, 
\[ \ts(gv)=\sigma(g)\ts(v)=g^{-1}v=-g^{-1}g^2v=-gv, \]
so
\[ \hat c^{\tilde\sigma}=\R(igv)\quad\mbox{and}
\quad\tth(igv)=-i\theta(g)\tth(v)=-igv.\]

We have shown that in the rank one case, 
associated to a noncompact imaginary root 
$\alpha$, a Cayley transformation $\mathbf c_\alpha = \tau(g)$ 
can be constructed so that it maps a given  
$\ts$- and $\tth$-stable Cartan subspace $c$ to a 
$\ts$- and $\tth$-stable Cartan subspace
$\hat c = \mathbf c_\alpha(c)$ such that 
the noncompact dimension of $\hat c^{\tilde\sigma}$
is one higher than that of $\cs$.
In the sequel, we want to generalize
this construction to an arbitrary polar orthogonal representation
$\tau:G\to O(V,\inn{\cdot}{\cdot})$. 

Indeed, suppose now that the rank of $\tau$ is arbitrary,
let $c$ be an arbitrary $\ts$- and $\tth$-stable Cartan subspace
and assume there exists a noncompact imaginary root 
$\alpha\in\mathscr A$ 
which we suppose fixed. Write $c=c_\alpha\oplus \C v_\alpha$ 
where $v_\alpha\in i(\cst)=c^{-\tilde\sigma}\cap c^{-\tilde\theta}$
is the co-root. Note that 
\[ \cs=c_\alpha^{\tilde\sigma}\oplus\R(iv_\alpha), \]
and $iv_\alpha\in c^{\tilde\theta}$. 
Now $(\Lg_\alpha,c\oplus\Lg_\alpha\cdot c)$ is a rank one polar 
action~\cite[Th.~2.12]{DK}. Since 
$V=c\oplus\bigoplus_{\alpha\in\mathscr A}\Lg_\alpha\cdot c$ is 
a $\inn{\cdot}{\cdot}$-orthogonal direct sum, $(\Lg_\alpha,c\oplus\Lg_\alpha\cdot c)$
is orthogonal with respect to the restriction of~$\inn{\cdot}{\cdot}$; 
we restrict
it to $(\Lg_\alpha,\C v_\alpha\oplus\tilde\Lg_\alpha\cdot v_\alpha)$ 
to get 
an irreducible polar orthogonal action of rank one. Since 
$X\in\Lg_\alpha\mapsto X\cdot v_\alpha$ is injective
on $\tilde\Lg_\alpha$, the kernel of 
this representation is contained in $\Lm$. Let $Z\in\Lm$. Then 
$Z\cdot v_\alpha=0$. If $Z\cdot\tilde\Lg_\alpha\cdot v_\alpha=0$, then
$[Z,\tilde\Lg_\alpha]\cdot v_\alpha=0$. Since 
$[Z,\tilde\Lg_\alpha]\subset\tilde\Lg_\alpha$, we get that 
$[Z,\tilde\Lg_\alpha]=0$, so $Z\in Z_{\mathfrak m}(\tilde\Lg_\alpha)$. 
Now $(\Lg'_\alpha,V_\alpha)$ is an effective
irreducible polar orthogonal action of rank one, where we have set
\[ \Lg'_\alpha=\Lg_\alpha/Z_{\mathfrak m}(\tilde\Lg_\alpha)\quad
\mbox{and}\quad V_\alpha = \C v_\alpha\oplus\tilde\Lg_\alpha\cdot v_\alpha. \]
Note that $\alpha$ can also be considered as a root of 
$(\Lg'_\alpha,V_\alpha)$, and then it is a noncompact imaginary 
root, so by the previous discussion we can find $g\in G_\alpha$ as above
and perform a Cayley transform $\mathbf c_\alpha=\tau(g)$ 
as follows:
\[ \hat c = \mathbf c_\alpha(c)= c_\alpha \oplus \C(g v_\alpha). \]
Note that
\[ \hat c^{\tilde\sigma}=c_\alpha^{\tilde\sigma}\oplus 
\R(gv_\alpha), \]
and $gv_\alpha\in\cmt$, so the noncompact dimension 
of $\hat c^{\tilde\sigma}$ is one higher than that of 
$\cs$. In a completely
analogous way, one can define a Cayley transform 
that increases the compact dimension of $\cs$ by one 
by using a compact real root.

\subsection{Uniqueness of extremal Cartan subspaces}

The Cayley transform allows us to derive some important 
properties of extremal Cartan subspaces. 

\begin{thm}\label{polar}
We have that $(K_\mathbf R,\Vr\cap iW)$ 
(resp.~$(K_\mathbf R,\Vr\cap W)$) 
is a polar representation. The sections are given by 
$\csmt$ (resp.~$\cst$), 
where $\cs$ is a maximally noncompact (resp.~compact) Cartan subspace 
of $\Vr=V^{\tilde\sigma}$. 
\end{thm}

\Pf It suffices to treat the case
of $(K_\mathbf R,\Vr\cap iW)$.
Let $\cs$ be a maximally noncompact Cartan subspace.
Then there are no noncompact imaginary 
roots, for otherwise a Cayley transform could be performed 
increasing the noncompact dimension of $\cs$. 
We claim that there exists $v_2\in\csmt$ such that 
\[ \Lk_{\mathbf R}(v_2)\oplus \csmt = V^{\tilde\sigma}\cap V^{-\tilde\theta}
=\Vr\cap iW. \]
In order to prove this claim, we first remark that~\cite[Thm.2.12]{DK}
\[ \Lg_v = \Lm +\sum_{\alpha(v)=0}\tilde\Lg_\alpha \]
for $v\in c$, 
\[ \Lg_v = \left\{ \begin{array}{ll}
 \Lm + \sum_{\alpha\ \mathrm{imag}}\tilde\Lg_\alpha 
& \mbox{for generic $v\in c^{-\tilde\omega}$,} \\
  \Lm + \sum_{\alpha\ \mathrm{real}}\tilde\Lg_\alpha 
& \mbox{for generic $v\in c^{\tilde\omega}$,} 
                    \end{array} \right. \] 
and
\[ \Lg_{v_1} = \underbrace{\left(\Lm^\omega+\sum_{\alpha\ \mathrm{real}}
\tilde\Lg_\alpha^\omega\right)}_{\subset\mathfrak k} \oplus 
\underbrace{\left(\Lm^{-\omega}+\sum_{\alpha\ \mathrm{real}}
\tilde\Lg_\alpha^{-\omega}\right)}_{\subset\mathfrak p} 
\quad\mbox{for generic $v_1\in\cst$,} \]
\[ \Lg_{v_2} = \underbrace{\left(\Lm^\omega+\sum_{\alpha\ \mathrm{imag}}
\tilde\Lg_\alpha^\omega\right)}_{\subset\mathfrak k} \oplus 
\underbrace{\Lm^{-\omega}}_{\subset\mathfrak p} 
\quad\mbox{for generic $v_2\in\csmt$,} \]
where in the last line we have used 
the nonexistence of noncompact imaginary 
roots. Select generic $v_1\in\cst$, $v_2\in\csmt$ and 
set $v=v_1+v_2\in\cs$. For each $\alpha\in\mathscr A$, 
\[ \alpha(v)=\underbrace{\alpha(v_1)}_{\in i\mathbf R}+
\underbrace{\alpha(v_2)}_{\in\mathbf R}, \] 
where at least one of the two summands on the right 
hand-side in not zero
by the choice of $v_1$, $v_2$. 
This shows that~$v$ is regular for $(G,V)$. By polarity, 
$\Lg\cdot v\oplus c = V$. Taking real parts in $\Vr$ yields
\[ \Lgr(v)\oplus \cs = \Vr, \]
which is the same as
\[ (\Lk_{\mathbf R}(v_1)+\Lp_{\mathbf R}(v_2))\oplus
(\Lk_{\mathbf R}(v_2)+\Lp_{\mathbf R}(v_1))\oplus \cst\oplus\csmt
=\Vr\cap W\oplus\Vr\cap iW. \]
In particular, 
\[ (\Lk_{\mathbf R}(v_2)+\Lp_{\mathbf R}(v_1))\oplus\csmt
=\Vr\cap iW. \]
The claim will follow if we show that $\Lk_{\mathbf R}(v_2)\supset
\Lp_{\mathbf R}(v_1)$. This is to be a consequence 
of $\Lk\cdot v_2\supset\Lp\cdot v_1$, as 
$\Lk\cdot v_2$ and $\Lp\cdot v_1$  
are $\ts$-stable and 
$\Lk_{\mathbf R}(v_2)=(\Lk\cdot v_2)^{\tilde\sigma}$,
$\Lp_{\mathbf R}(v_1)=(\Lp\cdot v_1)^{\tilde\sigma}$.

Now $\Lp\cdot v_1$ 
is spanned by
\[ 
\underbrace{\tilde\Lg_{\alpha}^{-\omega}}_{=0}\cdot v_1
\quad\mbox{for $\alpha$ imaginary, and} \]
\begin{eqnarray*}
(X_\alpha-\omega X_\alpha)\cdot v_1 &=& 
X_\alpha\cdot v_1 -\tom(X_\alpha\cdot v_1) \\
&=&\underbrace{\alpha(v_1)}_{\neq0}
(1-\tom)(X_\alpha\cdot v_\alpha)\quad\mbox{for
$\alpha$ complex and $X_\alpha\in\tilde\Lg_\alpha$.}
\end{eqnarray*}
On the other hand, $\Lk(v_2)$ is spanned by
\[ 
\tilde\Lg_{\alpha}^{\omega}\cdot v_2
\quad\mbox{for $\alpha$ real, and} \]
\begin{eqnarray*}
(X_\alpha+\omega X_\alpha)\cdot v_2 &=& 
X_\alpha\cdot v_2 -\tom(X_\alpha\cdot v_2) \\
&=&\underbrace{\alpha(v_2)}_{\neq0}
(1-\tom)(X_\alpha\cdot v_\alpha)\quad\mbox{for
$\alpha$ complex and $X_\alpha\in\tilde\Lg_\alpha$.}
\end{eqnarray*}
This proves that $\Lp\cdot v_1\subset\Lk\cdot v_2$, and hence 
that $\Lk_{\mathbf R}(v_2)\oplus \csmt =\Vr\cap iW$. 
Since $\inn{\Lg\cdot c}c=0$, we get that 
$\csmt$ is the $\inn{\cdot}{\cdot}$-orthogonal complement of 
$K_{\mathbf R}(v_2)$ in $\Vr\cap iW$. Since $K_{\mathbf R}$ is compact
and $\inn{\cdot}{\cdot}$ is positive-definite on $\Vr\cap iW$, 
it easily follows that $\csmt$ meets all the $K_{\mathbf R}$-orbits 
in $\Vr\cap iW$. Again by $\inn{\Lg\cdot c}c=0$, one has that
$\csmt$ meets all the other $K_{\mathbf R}$-orbits orthogonally. 
This finishes the proof. \EPf

\begin{cor}\label{extremal}
Two maximally noncompact (resp.~compact) Cartan subspaces 
$c_1^{\tilde\sigma}$ and $c_2^{\tilde\sigma}$  
of $V^{\tilde\sigma}=\Vr$ are $K_{\mathbf R}$-conjugate. 
As a consequence, there exists a unique $\Gr$-conjugacy class
of maximally noncompact (resp.~compact) Cartan subspaces
of $\Vr$.
\end{cor}

\Pf Again, it suffices to treat the case of 
maximally noncompact Cartan subspaces. 
By Theorem~\ref{polar}, we may assume that 
\[ c_1^{\tilde\sigma}\cap c_1^{-\tilde\theta}=
c_2^{\tilde\sigma}\cap c_2^{-\tilde\theta}. \]
Take a generic point $v_2$ lying therein. Since 
$v_2\in\Vr\cap iW$, we have that 
$\Lu_{v_2}=(\Lk_{\mathbf R})_{v_2}+(i\Lp_{\mathbf R})_{v_2}$,
and this is a decomposition into the $\pm1$-eigenspaces
of $\sigma$ on $\Lu_{v_2}$, so
\[ (U_{v_2})^\circ = (K_{\mathbf R})_{v_2}\exp[(i\Lp_{\mathbf R})_{v_2}]. \]
Consider the slice of the 
polar action $(U,V^{-\tilde\theta})$ at $v_2$; 
it is also polar with the same sections:
\[ c_i^{-\tilde\theta} = c_i^{-\tilde\sigma}\cap c_i^{-\tilde\theta}
\oplus c_i^{\tilde\sigma}\cap c_i^{-\tilde\theta} \]
for $i=1$, $2$. Now $c_1^{-\tilde\theta}$ 
and~$c_2^{-\tilde\theta}$ 
must be conjugate by an element of $(U_{v_2})^\circ$.
Since $\exp[(i\Lp_{\mathbf R})_{v_2}]$ centralizes $c$ (for 
$(i\Lp_{\mathbf R})_{v_2}=\Lm^{-\sigma}\cap\Lm^{\theta}$),
they must indeed be conjugate by an element of $(\Kr)_{v_2}$
(which necessarily 
fixes $c_1^{\tilde\sigma}\cap c_1^{-\tilde\theta}=
c_2^{\tilde\sigma}\cap c_2^{-\tilde\theta}$ since 
this is a section  of $(\Kr,\Vr\cap iW)$ and 
$v_2$ is a regular point of that action). 
Hence, so are 
$c_1^{\tilde\sigma}\cap c_1^{\tilde\theta}=
i(c_1^{-\tilde\sigma}\cap c_1^{-\tilde\theta})$ and 
$c_2^{\tilde\sigma}\cap c_2^{\tilde\theta}=
i(c_2^{-\tilde\sigma}\cap c_2^{-\tilde\theta})$. 
\mbox{ }\hfill\EPf

 \providecommand{\bysame}{\leavevmode\hbox to3em{\hrulefill}\thinspace}
\providecommand{\MR}{\relax\ifhmode\unskip\space\fi MR }
\providecommand{\MRhref}[2]{%
  \href{http://www.ams.org/mathscinet-getitem?mr=#1}{#2}
}
\providecommand{\href}[2]{#2}

%\bibliographystyle{amsalpha}
%\bibliography{ref} 

\end{document}